\documentclass[letterpaper,final,twocolumn]{IEEEtran}  

\usepackage{amsmath,mathrsfs}
\usepackage{amssymb}
\usepackage{color}
\usepackage[all]{xy}
\usepackage{graphicx,xspace,bm}
\usepackage[caption = false]{subfig}
\usepackage{url}
\usepackage{setspace}
\usepackage{enumitem}

\newtheorem{theorem}{Theorem}[section]

\newtheorem{lemma}[theorem]{Lemma}

\newtheorem{remark}[theorem]{Remark}
\newtheorem{example}[theorem]{Example}
\newtheorem{corollary}[theorem]{Corollary}
\newtheorem{proposition}[theorem]{Proposition}

\newcommand{\BB}{\mathcal{B}}
\renewcommand{\SS}{\mathcal{S}}
\newcommand{\II}{\mathcal{I}}
\newcommand{\JJ}{\mathcal{J}}
\newcommand{\KK}{\mathcal{K}}
\newcommand{\LL}{\mathcal{L}}
\newcommand{\UU}{\mathcal{U}}

\newcommand{\DD}{\mathcal{D}}

\newcommand{\EE}{\mathcal{E}}
\newcommand{\MM}{\mathcal{M}}

\newcommand{\XX}{\mathcal{X}}
\newcommand{\YY}{\mathcal{Y}}
\newcommand{\ZZ}{\mathcal{Z}}

\newcommand{\intr}{\mathrm{int}}

\newcommand{\cl}{\mathrm{cl}}

\newcommand{\abs}[1]{\ensuremath{\left\lvert{#1}\right\rvert}}
\newcommand{\norm}[1]{\ensuremath{\| #1 \|}}

\newcommand{\Bnorm}[1]{\ensuremath{\Big \| #1 \Big \|}}

\newcommand{\real}{{\mathbb{R}}}

\newcommand{\realnonnegative}{{\mathbb{R}}_{\ge 0}}

\newcommand{\integersnonnegative}{\mathbb{Z}_{\geq 0}}

\newcommand{\naturalnumbers}{\mathbb{N}}
\newcommand{\eps}{\epsilon}

\newcommand{\until}[1]{\{1,\dots,#1\}}
\newcommand{\map}[3]{#1:#2 \rightarrow #3}

\newcommand{\setdef}[2]{\{#1 \; | \; #2\}}

\newcommand{\set}[1]{\{#1\}}

\newcommand{\Lie}{\mathcal{L}}

\newcommand{\gradient}{\nabla}
\newcommand{\Eq}[1]{\mathrm{Eq}(#1)}

\newcommand{\levelset}[2]{#1^{-1}({\leq #2})}

\newcommand{\sdl}[1]{\mbox{Saddle}(#1)}
\newcommand{\lm}{\lambda}

\newcommand{\xo}{x_{*}}
\newcommand{\yo}{y_{*}}
\newcommand{\zo}{z_{*}}
\newcommand{\lmo}{\lambda_{*}}

\newcommand{\xt}{\tilde{x}}
\newcommand{\yt}{\tilde{y}}

\newcommand{\tht}{t}
\newcommand{\SD}{X_{\text{sp}}}
\newcommand{\SDp}{X_{\text{sp}}^{\text{p}}}
\newcommand{\PSD}{X_{\text{p-sp}}}

\newcommand{\proj}{\mathrm{proj}}
\newcommand{\HO}{\overline{H}}
\newcommand{\oprocendsymbol}{\hbox{$\bullet$}}
\newcommand{\oprocend}{\relax\ifmmode\else\unskip\hfill\fi\oprocendsymbol}

\newcommand{\longthmtitle}[1]{\textsl{(#1):}}

\newcommand{\changes}[1]{{\color{blue} #1}}
\renewcommand{\changes}[1]{#1}

\renewcommand{\theenumi}{(\roman{enumi})}

\allowdisplaybreaks

\begin{document}

\title{The role of convexity in saddle-point dynamics: Lyapunov
  function and robustness\thanks{A preliminary version of this work
    appeared at the 2016 Allerton Conference on Communication, Control,
    and Computing, Monticello, Illinois
    as~\cite{AC-EM-SL-JC:16-allerton}.}}

\author{Ashish Cherukuri \quad Enrique Mallada \quad Steven Low \quad
  Jorge Cort\'{e}s \thanks{Ashish Cherukuri and Jorge~Cort\'{e}s are
    with the Department of Mechanical and Aerospace Engineering,
    University of California, San Diego,
    \texttt{\{acheruku,cortes\}@ucsd.edu}, Enrique Mallada is with the
    Department of Electrical and Computer Engineering, Johns Hopkins
    University \texttt{mallada@hopkins.edu}, and Steven Low is with
    the Computing and Mathematical Sciences and the Electrical
    Engineering Department, California
    Institute of Technology, \texttt{slow@caltech.edu}}
}

\maketitle

\begin{abstract}
  This paper studies the projected saddle-point dynamics associated to
  a convex-concave function, which we term saddle function.  The
  dynamics consists of gradient descent of the saddle function in
  variables corresponding to convexity and (projected) gradient ascent
  in variables corresponding to concavity.  \changes{We examine the
    role that the local and/or global nature of the
    convexity-concavity properties of the saddle function plays in
    guaranteeing convergence and robustness of the dynamics. } Under
  the assumption that the saddle function is twice continuously
  differentiable, we provide a novel characterization of the
  omega-limit set of the trajectories of this dynamics in terms of the
  diagonal blocks of the Hessian. Using this characterization, we
  establish global asymptotic convergence of the dynamics under local
  strong convexity-concavity of the saddle function. When strong
  convexity-concavity holds globally, we establish three results.
  First, we identify a Lyapunov function \changes{(that decreases strictly along
  the trajectory)} for the projected
  saddle-point dynamics when the saddle function corresponds to the
  Lagrangian of a general constrained convex optimization
  problem. Second, \changes{for the particular case} when the saddle
  function is the Lagrangian of an \changes{equality-constrained}
  optimization problem, we show input-to-state stability of the
  saddle-point dynamics by providing an ISS Lyapunov function. Third,
  we \changes{use the latter result to} design an opportunistic
  state-triggered implementation of the dynamics. Various examples
  illustrate our results.
\end{abstract}

\section{Introduction}

Saddle-point dynamics and its variations have been used extensively in
the design and analysis of distributed feedback controllers and
optimization algorithms in several domains, including power networks,
network flow problems, and zero-sum games.  The analysis of the global
convergence of this class of dynamics typically relies on some global
strong/strict convexity-concavity property of the saddle function
defining the dynamics. The main aim of this paper is to refine this
analysis by unveiling two ways in which convexity-concavity of the
saddle function plays a role. First, we show that local strong
convexity-concavity is enough to conclude global asymptotic
convergence, thus generalizing previous results that rely on global
strong/strict convexity-concavity instead.  \changes{ Second, we show
  that, if global strong convexity-concavity holds, then one can
  identify a novel Lyapunov function for the projected saddle-point
  dynamics for the case when the saddle function is the Lagrangian of
  a constrained optimization problem. This, in turn, implies a
  stronger form of convergence, that is, input-to-state stability
  (ISS) and has important implications in the practical implementation
  of the saddle-point dynamics.}

\subsubsection*{Literature review}

The analysis of the convergence properties of (projected) saddle-point
dynamics to the set of saddle points goes back
to~\cite{TK:56,KA-LH-HU:58}, motivated by the study of nonlinear
programming and optimization. These works employed direct methods,
examining the approximate evolution of the distance of the
trajectories to the saddle point and concluding attractivity by
showing it to be decreasing. Subsequently, motivated by the extensive
use of the saddle-point dynamics in congestion control problems, the
literature on communication networks developed a Lyapunov-based and
passivity-based asymptotic stability analysis, see
e.g.~\cite{MC-SHL-ARC-JCD:07} and references therein.  Motivated by
network optimization, more recent works~\cite{DF-FP:10,JW-NE:11} have
employed indirect, LaSalle-type arguments to analyze asymptotic
convergence.  For this class of problems, the aggregate nature of the
objective function and the local computability of the constraints make
the saddle-point dynamics corresponding to the Lagrangian naturally
distributed.  Many other works exploit this dynamics to solve network
optimization problems for various applications, e.g., distributed
convex optimization~\cite{JW-NE:11,BG-JC:14-tac}, distributed linear
programming~\cite{DR-JC:15-tac}, bargaining
problems~\cite{DR-JC:16-tcns}, and power
networks~\cite{XZ-AP:15-auto,CZ-UT-NL-SL:14,NL-CZ-LC:16,TS-CDP-AJS:17-tac,EM-CZ-SL:14-arxiv}.
Another area of application is game theory, where saddle-point
dynamics is applied to find the Nash equilibria of two-person zero-sum
games~\cite{BG-JC:13-auto,LJR-SAB-SSS:16-tac}. In the context of
distributed optimization, the recent work~\cite{SKN-JC:16-sicon}
employs a (strict) Lyapunov function approach to ensure asymptotic
convergence of saddle-point-like dynamics. The work~\cite{TH-IL:14}
examines the asymptotic behavior of the saddle-point dynamics when the
set of saddle points is not asymptotically stable and, instead,
trajectories exhibit oscillatory behavior.  Our previous work has
established global asymptotic convergence of the saddle-point
dynamics~\cite{AC-BG-JC:17-sicon} and the projected saddle-point
dynamics~\cite{AC-EM-JC:16-scl} under global strict
convexity-concavity assumptions. The works mentioned above require
similar or stronger global assumptions on the convexity-concavity
properties of the saddle function to ensure convergence.  Our results
here directly generalize the convergence properties reported above.
Specifically, we show that traditional assumptions on the problem
setup can be relaxed if convergence of the dynamics is the desired
property: global convergence of the projected saddle-point dynamics
can be guaranteed under local strong convexity-concavity assumptions.
Furthermore, if traditional assumptions do hold, then a stronger
notion of convergence, that also implies robustness, is guaranteed: if
strong convexity-concavity holds globally, the dynamics admits a
Lyapunov function and in the absence of projection, the dynamics is
ISS, admitting an ISS Lyapunov function.

\subsubsection*{Statement of contributions}

Our starting point is the definition of the projected saddle-point
dynamics for a differentiable convex-concave function, referred to as
saddle function. The dynamics has three components: gradient descent,
projected gradient ascent, and gradient ascent of the saddle function,
where each gradient is with respect to a subset of the arguments of
the function. This unified formulation encompasses all forms of the
saddle-point dynamics mentioned in the literature review above.  Our
contributions shed light on the effect that the convexity-concavity of
the saddle function has on the convergence attributes of the projected
saddle-point dynamics.  Our first contribution is a novel
characterization of the omega-limit set of the trajectories of the
projected saddle-point dynamics in terms of the diagonal Hessian
blocks of the saddle function.  To this end, we use the distance to a
saddle point as a LaSalle function, express the Lie derivative of this
function in terms of the Hessian blocks, and show it is nonpositive
using second-order properties of the saddle function.  Building on
this characterization, our second contribution establishes global
asymptotic convergence of the projected saddle-point dynamics to a
saddle point assuming only local strong convexity-concavity of the
saddle function.  Our third contribution identifies a novel Lyapunov
function for the projected saddle-point dynamics for the case when
strong convexity-concavity holds globally and the saddle function can
be written as the Lagrangian of a constrained optimization problem.
This discontinuous Lyapunov function can be interpreted as multiple
continuously differentiable Lyapunov functions, one for each set in a
particular partition of the domain determined by the projection
operator of the dynamics. Interestingly, the identified Lyapunov
function is the sum of two previously known and independently
considered LaSalle functions.  When the saddle function takes the form
of the Lagrangian of an equality constrained optimization, then no
projection is present.  In such scenarios, if the saddle function
satisfies global strong convexity-concavity, our fourth contribution
establishes input-to-state stability (ISS) of the dynamics with
respect to the saddle point by providing an ISS Lyapunov function.
Our last contribution uses this function to design an opportunistic
state-triggered implementation of the saddle-point dynamics. We show
that the trajectories of this discrete-time system converge
asymptotically to the saddle points and that executions are Zeno-free,
i.e., that the difference between any two consecutive triggering times
is lower bounded by a common positive quantity.  Examples illustrate
our results.

\section{Preliminaries}\label{sec:prelims}
This section introduces our notation and preliminary notions on
convex-concave functions, discontinuous dynamical systems, and
input-to-state stability.

\subsection{Notation}
Let $\real$, $\realnonnegative$, and $\naturalnumbers$ denote the set
of real, nonnegative real, and natural numbers, respectively. We let
$\norm{\cdot}$ denote the $2$-norm on $\real^n$ and the respective
induced norm on $\real^{n \times m}$.  Given $x,y\in \real^n$, $x_i$
denotes the $i$-th component of $x$, and $x \le y$ denotes $x_i \le
y_i$ for $i \in \until{n}$.  For vectors $u \in \real^n$ and $w \in
\real^m$, the vector $(u;w) \in \real^{n+m}$ denotes their
concatenation.  For $a \in \real$ and $b \in \realnonnegative$, we let
\begin{align*}
  [a]_{b}^+ = \begin{cases} a, & \quad \text{ if } b> 0,
    \\
    \max\{0,a\}, & \quad \text{ if } b = 0.
  \end{cases}
\end{align*}
For vectors $a \in \real^n$ and $b \in \realnonnegative^n$, $[a]_{b}^+$
denotes the vector whose $i$-th component is $[a_i]_{b_i}^+$, for $i \in
\until{n}$.  Given a set $\SS \subset \real^n$, we denote by $\cl(\SS)$,
$\intr(\SS)$, and $\abs{\SS}$ its closure, interior, and cardinality,
respectively.  The distance of a point $x \in \real^n$ to the set $\SS
\subset \real^n$ in $2$-norm is $\norm{x}_\SS = \inf_{y \in \SS}
\norm{x-y}$.  The projection of $x$ onto a closed set $\SS$ is defined
as the set $\proj_{\SS}(x) = \setdef{y \in \SS}{\norm{x-y} =
\norm{x}_{\SS}}$. When $\SS$ is also convex, $\proj_{\SS}(x)$ is a singleton
for any $x \in \real^n$.  For a matrix $A \in \real^{n \times n}$, we
use $A \succeq 0$, $A \succ 0$, $A \preceq 0$, and $A \prec 0$ to denote
that $A$ is positive semidefinite, positive definite, negative
semidefinite, and negative definite, respectively. For a symmetric
matrix $A \in \real^{n \times n}$, $\lm_{\min}(A)$ and $\lm_{\max}(A)$
denote the minimum and maximum eigenvalue of $A$. For a real-valued
function $F:\real^n \times \real^m \to \real$, $(x,y) \mapsto F(x,y)$,
we denote by $\gradient_x F$ and $\gradient_y F$ the column vector of
partial derivatives of $F$ with respect to the first and second
arguments, respectively.  Higher-order derivatives follow the convention
$\gradient_{xy} F = \frac{\partial^2 F}{\partial x \partial y}$,
$\gradient_{xx} F = \frac{\partial^{2} F}{\partial x^2}$, and so on.  A
function $\map{\alpha}{\realnonnegative}{\realnonnegative}$ is class
$\KK$ if it is continuous, strictly increasing, and $\alpha(0) = 0$.
The set of unbounded class $\KK$ functions are called $\KK_\infty$
functions.  A function $\map{\beta}{\realnonnegative \times
\realnonnegative}{\realnonnegative}$ is class $\KK \LL$ if for any $t
\in \realnonnegative$, $x \mapsto \beta(x,t)$ is class $\KK$ and for any
$x \in \realnonnegative$, $t \mapsto \beta(x,t)$ is continuous,
decreasing with $\beta(t,x) \to 0$ as $t \to \infty$.

\subsection{Saddle points and convex-concave
functions}\label{subsec:saddle-points}

Here, we review notions of convexity, concavity, and saddle points
from~\cite{SB-LV:09}. A function $\map{f}{\XX}{\real}$ is
\emph{convex} if
\begin{align*}
  f(\lm x + (1-\lm)x') \le \lm f(x) + (1-\lm) f(x'),
\end{align*}
for all $x,x' \in \XX$ (where $\XX$ is a convex domain) 
and all $\lm \in [0,1]$. A convex
differentiable $f$ satisfies the following \emph{first-order
  convexity condition}
\begin{align*}
  f(x') \ge f(x) + (x' - x)^\top \gradient f(x),
\end{align*}
for all $x, x' \in \XX$.  A twice differentiable function $f$ is
\emph{locally strongly convex} at $x \in \XX$ if $f$ is convex and
$\gradient^2 f(x) \succeq mI$ for some $m > 0$ \changes{(note that
  this is equivalent to having $\gradient^2 f \succ 0$ in a
  neighborhood of~$x$).} Moreover, a twice differentiable $f$ is
\emph{strongly convex} if $\gradient^2 f(x) \succeq mI$ for all $x \in
\XX$ for some $m > 0$. A function $\map{f}{\XX}{\real}$ is
\emph{concave}, \emph{locally strongly concave}, or \emph{strongly
  concave} if $-f$ is convex, locally strongly convex, or strongly
convex, respectively.  A function $\map{F}{\XX \times \YY}{\real}$ is
\emph{convex-concave} (on $\XX \times \YY$) if, given any point
$(\xt,\yt) \in \XX \times \YY$, $x \mapsto F(x,\yt)$ is convex and $y
\mapsto F(\xt,y)$ is concave. When the space $\XX \times \YY$ is clear
from the context, we refer to this property as $F$ being
convex-concave in $(x,y)$. A point $(\xo,\yo) \in \XX \times \YY$ is a
\emph{saddle point} of $F$ on the set $\XX \times \YY$ if $F(\xo,y)
\le F(\xo,\yo) \le F(x,\yo)$, for all $x \in \XX$ and $y \in \YY$. The
set of saddle points of a convex-concave function~$F$ is convex. The
function $F$ is \emph{locally strongly convex-concave} at a saddle
point $(x,y)$ if it is convex-concave and either $\gradient_{xx}
F(x,y) \succeq mI$ or $\gradient_{yy} F(x,y) \preceq -mI$ for some $m
> 0$. Finally, $F$ is \emph{globally strongly convex-concave} if it is
convex-concave and either $x \mapsto F(x,y)$ is strongly convex for
all $y \in \YY$ or $y \mapsto F(x,y)$ is strongly concave for all $x
\in \XX$.

\subsection{Discontinuous dynamical systems}\label{subsec:disc}

Here we present notions of discontinuous dynamical
systems~\cite{AB-FC:06,JC:08-csm-yo}. Let $\map{f}{\real^n}{\real^n}$ be
Lebesgue measurable and locally bounded.  Consider the differential
equation
\begin{equation}\label{eq:dis-dyn}
  \dot x = f(x) .
\end{equation}
A map $\map{\gamma}{[0,T)}{\real^n}$ is a \emph{(Caratheodory)
  solution} of~\eqref{eq:dis-dyn} on the interval $[0,T)$ if it is
absolutely continuous on $[0,T)$ and satisfies $\dot \gamma(t) =
f(\gamma(t))$ almost everywhere in $[0,T)$.  We use the terms solution
and trajectory interchangeably.  A set $\SS \subset \real^n$ is
\emph{invariant} under~\eqref{eq:dis-dyn} if every solution starting
in $\SS$ remains in $\SS$.  For a solution $\gamma$
of~\eqref{eq:dis-dyn} defined on the time interval $[0,\infty)$, the
\emph{omega-limit} set $\Omega(\gamma)$ is defined by
\begin{multline*}
  \Omega(\gamma) = \setdef{y \in \real^n}{\exists
    \{t_k\}_{k=1}^{\infty} \subset [0,\infty) \text{ with } \lim_{k
      \to \infty} t_k = \infty
    \\
    \text{ and } \lim_{k \to \infty} \gamma(t_k) = y} .
\end{multline*}
If the solution $\gamma$ is bounded, then $\Omega (\gamma) \neq
\emptyset$ by the Bolzano-Weierstrass theorem~\cite[p. 33]{SL:93}.
Given a continuously differentiable function
$\map{V}{\real^n}{\real}$, the \emph{Lie derivative of $V$
  along~\eqref{eq:dis-dyn}} at $x \in \real^n$ is $\Lie_f V(x) =
\gradient V(x)^\top f(x)$.  The next result is a simplified version
of~\cite[Proposition 3]{AB-FC:06}.

\begin{proposition}\longthmtitle{Invariance principle for
    discontinuous Caratheodory systems}\label{pr:invariance-cara}
  Let $\SS \in \real^n$ be compact and invariant. Assume that, for
  each point $x_0 \in \SS$, there exists a unique solution
  of~\eqref{eq:dis-dyn} starting at $x_0$ and that its omega-limit set
  is invariant too. Let $\map{V}{\real^n}{\real}$ be a continuously
  differentiable map such that $\Lie_f V(x) \le 0$ for all $x \in
  \SS$. Then, any solution of~\eqref{eq:dis-dyn} starting at $\SS$
  converges to the largest invariant set in $\cl(\setdef{x \in
    \SS}{\Lie_f V(x) = 0})$.
\end{proposition}

\subsection{Input-to-state stability}\label{subsec:iss-prelim}

Here, we review the notion of input-to-state stability (ISS)
following~\cite{YL-ES-YW:95}. Consider a system
\begin{equation}\label{eq:iss-sys}
  \dot x = f(x,u),
\end{equation}
where $x \in \real^n$ is the state,
$\map{u}{\realnonnegative}{\real^m}$ is the input that is measurable
and locally essentially bounded, and $\map{f}{\real^n \times
  \real^m}{\real^n}$ is locally Lipschitz. Assume that starting from
any point in $\real^n$, the trajectory of~\eqref{eq:iss-sys} is
defined on $\realnonnegative$ for any given control. Let $\Eq{f}
\subset \real^n$ be the set of equilibrium points of the unforced
system. Then, the system~\eqref{eq:iss-sys} is \emph{input-to-state
  stable} (ISS) with respect to $\Eq{f}$ if there exists $\beta \in
\KK \LL$ and $\gamma \in \KK$ such that each trajectory $t \mapsto
x(t)$ of~\eqref{eq:iss-sys} satisfies
\begin{align*}
  \norm{x(t)}_{\Eq{f}} \le \beta(\norm{(x(0)}_{\Eq{f}},t) 
   + \gamma(\norm{u}_{\infty})
\end{align*}
for all $t \ge 0$, where $\norm{u}_\infty = \text{ess sup}_{t \ge 0}
\norm{u(t)}$ is the essential supremum (see~\cite[p. 185]{SL:93} for
the definition) of~$u$. This notion captures the graceful degradation
of the asymptotic convergence properties of the unforced system as the
size of the disturbance input grows.  One convenient way of showing
ISS is by finding an ISS-Lyapunov function. An
\emph{ISS-Lyapunov function} with respect to the set $\Eq{f}$ for
system~\eqref{eq:iss-sys} is a differentiable function
$\map{V}{\real^n}{\realnonnegative}$ such that
\begin{enumerate}
  \item there exist $\alpha_1, \alpha_2 \in \KK_\infty$ such that for
    all $x \in \real^n$, 
    \begin{align}\label{eq:iss-lyap-c1}
      \alpha_1(\norm{x}_{\Eq{f}}) \le V(x) \le
      \alpha_2(\norm{x}_{\Eq{f}});
    \end{align}
  \item there exists a continuous, positive definite function
    $\map{\alpha_3}{\realnonnegative}{\realnonnegative}$ and $\gamma
    \in \KK_\infty$ such that
    \begin{align}\label{eq:iss-lyap-c2}
      \gradient V(x)^\top f(x,v) \le - \alpha_3(\norm{x}_{\Eq{f}})
    \end{align}
    for all $x \in \real^n$, $v \in \real^m$ for which
    $\norm{x}_{\Eq{f}} \ge \gamma(\norm{v})$.
\end{enumerate}

\begin{proposition}\longthmtitle{ISS-Lyapunov function
  implies ISS}\label{pr:iss-lyap-iss}
  If~\eqref{eq:iss-sys} admits an ISS-Lyapunov function, then it is ISS.
\end{proposition}

\section{Problem statement}\label{sec:problem}

In this section, we provide a formal statement of the problem of
interest.  Consider a twice continuously differentiable function
$\map{F}{\real^n \times \realnonnegative^p \times \real^m}{\real}$,
$(x,y,z) \mapsto F(x,y,z)$, which we refer to as \emph{saddle
  function}. With the notation of Section~\ref{subsec:saddle-points},
we set $\XX = \real^n$ and $\YY = \realnonnegative^p \times
\real^m$, and assume that $F$ is convex-concave on $(\real^n) \times
(\realnonnegative^p \times \real^m)$. Let $\sdl{F}$ denote its
(non-empty) set of
saddle points. We define the
\emph{projected saddle-point dynamics} for $F$ as
\begin{subequations}\label{eq:proj-saddle-dyn}
  \begin{align}
    \dot x &= -\gradient_x F(x,y,z), \label{eq:proj-saddle-dyn-1}
    \\
    \dot y &= [\gradient_y F(x,y,z)]_y^+, \label{eq:proj-saddle-dyn-2}
    \\
    \dot z &= \gradient_z F(x,y,z). \label{eq:proj-saddle-dyn-3}
  \end{align}
\end{subequations}
When convenient, we use the map $\map{\PSD}{\real^n \times
  \realnonnegative^p \times \real^m}{\real^n \times \real^p \times
  \real^m}$ to refer to the dynamics~\eqref{eq:proj-saddle-dyn}.  Note
that the domain $\real^n \times \realnonnegative^p \times \real^m$ is
invariant under $\PSD$ (this follows from the definition of the
projection operator) and its set of equilibrium points precisely
corresponds to $ \sdl{F}$ (this follows from the defining property of
saddle points and the first-order condition for convexity-concavity of
$F$). Thus, a saddle point $(\xo,\yo,\zo)$ satisfies
\begin{subequations}\label{eq:equi-cond}
  \begin{align}
    \gradient_x F(\xo,\yo,\zo) = 0, \quad \gradient_z F(\xo,\yo,\zo) =
    0, 
    \\
    \gradient_y F(\xo,\yo,\zo) \le 0, \quad \yo^\top \gradient_y
    F(\xo,\yo,\zo) = 0.
  \end{align}
\end{subequations}
Our interest in the dynamics~\eqref{eq:proj-saddle-dyn} is motivated
by two bodies of work in the literature: one that analyzes primal-dual
dynamics, corresponding to~\eqref{eq:proj-saddle-dyn-1} together
with~\eqref{eq:proj-saddle-dyn-2}, for solving inequality constrained
network optimization problems, see
e.g.,~\cite{KA-LH-HU:58,DF-FP:10,EM-CZ-SL:14-arxiv,CZ-UT-NL-SL:14};
and the other one analyzing saddle-point dynamics, corresponding
to~\eqref{eq:proj-saddle-dyn-1} together
with~\eqref{eq:proj-saddle-dyn-3}, for solving equality constrained
problems and finding Nash equilibrium of zero-sum games, see
e.g.,~\cite{AC-BG-JC:17-sicon} and references therein.  By
considering~\eqref{eq:proj-saddle-dyn-1}-\eqref{eq:proj-saddle-dyn-3}
together, we aim to unify these lines of work.
\changes{Below we explain further the significance of the dynamics in
solving specific network optimization problems.

\begin{remark}\longthmtitle{Motivating examples}\label{re:motivation}
  {\rm Consider the following constrained convex optimization problem
    \begin{align*}
      \min \setdef{ f(x)}{ g(x) \le 0, \, Ax = b},
    \end{align*}
    where $\map{f}{\real^n}{\real}$ and $\map{g}{\real^n}{\real^p}$
    are convex continuously differentiable functions, $A \in \real^{m
      \times n}$, and $b \in \real^m$. Under zero duality gap, saddle
    points of the associated Lagrangian $L(x,y,z) = f(x) + y^\top g(x)
    + z^\top (Ax - b)$ correspond to the primal-dual optimizers of the
    problem.  This observation motivates the search for the saddle
    points of the Lagrangian, which can be done via the projected
    saddle-point dynamics~\eqref{eq:proj-saddle-dyn}.  In many network
    optimization problems, $f$ is the summation of individual costs of
    agents and the constraints, defined by $g$ and $A$, are such that
    each of its components is computable by one agent interacting with
    its neighbors. This structure renders the projected saddle-point
    dynamics of the Lagrangian implementable in a distributed
    manner. Motivated by this, the dynamics is widespread in network
    optimization scenarios. For example, in optimal dispatch of power
    generators~\cite{CZ-UT-NL-SL:14,NL-CZ-LC:16,TS-CDP-AJS:17-tac,EM-CZ-SL:14-arxiv},
    the objective function is the sum of the individual cost function
    of each generator, the inequalities consist of generator capacity
    constraints and line limits, and the equality encodes the power
    balance at each bus.  In congestion control of communication
    networks~\cite{MC-SHL-ARC-JCD:07, FP-EM:09, DF-FP:10}, the cost function is
    the summation of the negative of the utility of the communicated
    data, the inequalities define constraints on channel capacities,
    and equalities encode the data balance at each node.  \oprocend }
\end{remark}
}

Our main objectives are to identify conditions that guarantee that the
set of saddle points is globally asymptotically stable under the
dynamics~\eqref{eq:proj-saddle-dyn} and formally characterize the
robustness properties using the concept of input-to-state stability.
\changes{ The rest of the paper is structured as
  follows. Section~\ref{sec:proj-saddle} investigates novel conditions
  that guarantee global asymptotic convergence relying on LaSalle-type
  arguments. Section~\ref{sec:Lyapunov-function} instead identifies a
  strict Lyapunov function for constrained convex optimization
  problems. This finding allows us in Section~\ref{sec:ISS+self} to go
  beyond convergence guarantees and explore the robustness properties
  of the saddle-point dynamics.}

\section{Local properties of the saddle function imply global
  convergence}\label{sec:proj-saddle}

Our first result of this section provides a novel characterization of
the omega-limit set of the trajectories of the projected saddle-point
dynamics~\eqref{eq:proj-saddle-dyn}.

\begin{proposition}\longthmtitle{Characterization of the omega-limit
    set of solutions of $\PSD$}\label{prop:omega-charac-proj}
  Given a twice continuously differentiable, convex-concave function
  $F$, \changes{each point in the set $\sdl{F}$} is stable under the projected saddle-point
  dynamics $\PSD$ and the omega-limit set of every solution is contained
  in the largest invariant set $\MM$ in $\EE(F)$, where
  \begin{align}\label{eq:E-proj}
    \EE(F) & = \setdef{(x,y,z) \in \real^n \times \realnonnegative^p
      \times \real^m}{ \notag
      \\
      & \quad (x -\xo; y- \yo; z-\zo) \in
      \ker(\HO(x,y,z,\xo,\yo,\zo)), \notag
      \\
      & \quad \text{ for all } (\xo,\yo,\zo) \in \sdl{F} },
  \end{align}
  and
  \begin{flalign}\label{eq:Hdefs}
    & \HO(x,y,z,\xo,\yo,\zo) = \int_0^1  H(x(s),y(s),z(s)) ds,  \notag
    \\
    & (x(s),y(s),z(s)) = (\xo,\yo,\zo) + s(x-\xo,y-\yo,z-\zo), \notag
    \\
    & H(x,y,z) = \left[\begin{array}{ccc}
        -\gradient_{xx} F & 0 & 0\\
        0 & \gradient_{yy} F & \gradient_{yz} F \\
        0 & \gradient_{zy} F & \gradient_{zz} F
      \end{array}\right]_{(x,y,z)}.
  \end{flalign}
\end{proposition}
\begin{IEEEproof}
  The proof follows from the application of the LaSalle Invariance
  Principle for discontinuous Caratheodory systems
  (cf. Proposition~\ref{pr:invariance-cara}). Let $(\xo,\yo,\zo) \in
  \sdl{F}$ and $\map{V_1}{\real^n \times \realnonnegative^p \times
    \real^m}{\realnonnegative}$ be defined as
  \begin{align}\label{eq:V1}
    V_1(x,y,z) \!=\! \frac{1}{2} \bigl( \norm{x-\xo}^2 \!+\! \norm{y - \yo}^2
    \!+\! \norm{z-\zo}^2 \bigr).
  \end{align}
  The Lie derivative of $V_1$ along~\eqref{eq:proj-saddle-dyn} is
  \begin{align}
    & \Lie_{\PSD} V_1(x,y,z) \notag
    \\
    & \quad = - (x-\xo)^\top \gradient_x F(x,y,z) + (y-\yo)^\top
    [\gradient_y F(x,y,z)]_y^+ \notag
    \\
    & \qquad + (z - \zo)^\top \gradient_z F(x,y,z) \notag
    \\
    & \quad = - (x-\xo)^\top \gradient_x F(x,y,z) + (y-\yo)^\top
    \gradient_y F(x,y,z) \notag
    \\
    & \qquad + (z - \zo)^\top \gradient_z F(x,y,z)\notag
    \\
    & \qquad + (y - \yo)^\top ( [\gradient_y F(x,y,z)]_y^+ -
    \gradient_y F(x,y,z))\notag
    \\
    & \quad \le - (x-\xo)^\top \gradient_x F(x,y,z) + (y-\yo)^\top
    \gradient_y F(x,y,z)\notag
    \\
    & \qquad + (z - \zo)^\top \gradient_z F(x,y,z) ,
    \label{eq:lie-v1-ineq}
  \end{align}
  where the last inequality follows from the fact that $T_i = (y -
  \yo)_i([\gradient_y F(x,y,z)]_y^+ - \gradient_y F(x,y,z))_i \le 0$
  for each $i \in \until{p}$. Indeed if $y_i > 0$, then $T_i = 0$ and
  if $y_i = 0$, then $(y - \yo)_i \le 0$ and $([\gradient_y
  F(x,y,z)]_y^+ - \gradient_y F(x,y,z))_i \ge 0$ which implies that
  $T_i \le 0$. Next, denoting $\lm = (y;z)$ and $\lmo = (\yo,\zo)$, we
  simplify the above inequality as
  \begin{align*}
    & \Lie_{\PSD} V_1(x,y,z)
    \\
    & \quad \le - (x-\xo)^\top \gradient_x F(x,\lm) 
    + (\lm - \lmo)^\top \gradient_{\lm} F(x,\lm) 
    \\
    & \quad \overset{(a)}{=}-(x-\xo)^\top \int_0^1 
    \Bigl( \gradient_{xx} F(x(s),\lm(s)) (x-\xo)  
    \\
    & \qquad + \gradient_{\lm x} F(x(s),\lm(s)) 
    (\lm-\lmo) \Bigr) ds 
    \\
    &\qquad+ (\lm - \lmo)^\top \int_0^1 \Bigl( \gradient_{x \lm} 
    F(x(s),\lm(s)) (x-\xo)  
    \\
    & \qquad + \gradient_{\lm \lm } F(x(s),\lm(s)) 
    (\lm - \lmo) \Bigr) ds 
    \\
    &\quad \overset{(b)}{=}[x-\xo;\lm-\lmo]^\top \HO(x,\lm,\xo,\lmo)
    \left[\begin{array}{c} x-\xo
        \\
        \lm-\lmo \end{array} \right]
    \overset{(c)}{\leq} 0,
  \end{align*}
  where (a) follows from the fundamental theorem of calculus using the
  notation $x(s) = \xo + s(x - \xo)$ and $\lm(s) = \lmo + s(\lm -
  \lmo)$ and recalling from~\eqref{eq:equi-cond} that $\gradient_x
  F(\xo,\lmo) = 0$ and $(\lm - \lmo)^\top \gradient_\lambda
  F(\xo,\lmo) \le 0$; (b) follows from the definition of $\HO$ using
  $(\gradient_{\lm x} F(x,\lm))^\top = \gradient_{x \lm} F(x,\lm)$;
  and (c) follows from the fact that $\HO$ is negative
  semi-definite. Now using this fact that $\Lie_{\PSD}V_1$ is
  nonpositive at any point, one can deduce, see e.g.~\cite[Lemma
  4.2-4.4]{AC-EM-JC:16-scl}, that starting from any point
  $(x(0),y(0),z(0))$ a unique trajectory of $\PSD$ exists, is
  contained in the compact set $V_1^{-1}(V_1(x(0),y(0),z(0))) \cap
  (\real^n \times \realnonnegative^p \times \real^m)$ at all times,
  and its omega-limit set is invariant. These facts imply that the
  hypotheses of Proposition~\ref{pr:invariance-cara} hold and so, we
  deduce that the solutions of the dynamics $\PSD$ converge to the
  largest invariant set where the Lie derivative is zero, that is, the
  set
  \begin{align}\label{eq:E-proj-xyz}
    & \EE(F,\xo, \yo, \zo) = \setdef{(x,y,z) \in \real^n \times
      \realnonnegative^p \times \real^m}{ \notag \\ & (x;y;z)
	-(\xo;\yo;\zo) \in \ker(\HO(x,y,z,\xo,\yo,\zo))}.
  \end{align}
  Finally, since $(\xo,\yo,\zo)$ was chosen arbitrary, we get that the
  solutions converge to the largest invariant set $\MM$ contained in
  $\EE(F) = \bigcap_{(\xo,\yo,\zo)\in\sdl{F}} \EE(F,\xo,\yo,\zo)$,
  concluding the proof. 
\end{IEEEproof}

Note that the proof of Proposition~\ref{prop:omega-charac-proj} shows
that the Lie derivative of the function $V_1$ is negative, but not
strictly negative, outside the set~$\sdl{F}$.  
\changes{From Proposition~\ref{prop:omega-charac-proj} and the
definition~\eqref{eq:E-proj}, we deduce that if a point 
$(x,y,z)$ belongs to the omega-limit set
(and is not a saddle point), then 
the line integral of the Hessian block
matrix~\eqref{eq:Hdefs} from the any saddle point to $(x,y,z)$
cannot be full rank. Elaborating further,    
\begin{enumerate}
  \item if $\gradient_{xx} F$ is full rank at a saddle point
    $(x_*,y_*,z_*)$ and if the point $(x,y,z) \not \in \sdl{F}$ 
    belongs to the omega-limit set, then $x = x_*$, and  
  \item if $\begin{bmatrix} \gradient_{yy} F & \gradient_{yz} F \\
      \gradient_{zy} F & \gradient_{zz} F \end{bmatrix}$  is full rank
    at a saddle point $(x_*,y_*,z_*)$, then $(y,z) = (y_*,z_*)$.
\end{enumerate}
These properties are used in the next result which shows that local
strong convexity-concavity at a saddle point together
with global convexity-concavity of the saddle function are enough to
guarantee global convergence.proving Theorem 4.2.}

\begin{theorem}\longthmtitle{Global asymptotic stability of the set of
    saddle points under~$\PSD$}\label{th:proj-local2global}
  Given a twice continuously differentiable, convex-concave function $F$
  which is locally strongly convex-concave at a saddle point, the set
  $\sdl{F}$ is globally asymptotically stable under the projected
  saddle-point dynamics $\PSD$ and the convergence of trajectories is to
  a point.
\end{theorem}
\begin{IEEEproof}
  Our proof proceeds by characterizing the set $\EE(F)$ defined
  in~\eqref{eq:E-proj}. Let $(\xo,\yo,\zo)$ be a saddle point at which
  $F$ is locally strongly convex-concave.  Without loss of
  generality, assume that $\gradient_{xx} F(\xo,\yo,\zo) \succ 0$ (the
  case of negative definiteness of the other Hessian block can be
  reasoned analogously).
  Let $(x,y,z) \in \EE(F,\xo,\yo,\zo)$ (recall the definition of this
  set in~\eqref{eq:E-proj-xyz}).  Since $\gradient_{xx} F(\xo,\yo,\zo)
  \succ 0$ and $F$ is twice continuously differentiable, we have that
  $\gradient_{xx} F$ is positive definite in a neighborhood of
  $(\xo,\yo,\zo)$ and so
  \begin{align*}
    \int_0^1 \gradient_{xx} F(x(s),y(s),z(s)) ds \succ 0,
  \end{align*}
  where $x(s) = \xo+s(x-\xo)$, $y(s) = \yo+s(y-\yo)$, and $z(s) = \zo
  + s(z-\zo)$.  Therefore, by definition of $\EE(F,\xo,\yo,\zo)$, it
  follows that $x = \xo$ and so, $\EE(F,\xo,\yo,\zo) \subseteq \{\xo\}
  \times (\realnonnegative^p \times \real^m)$.  From
  Proposition~\ref{prop:omega-charac-proj} the trajectories of $\PSD$
  converge to the largest invariant set $\MM$ contained in
  $\EE(F,\xo,\yo,\zo)$. To characterize this set, let $(\xo,y,z) \in
  \MM$ and $t \mapsto (\xo,y(t),z(t))$ be a trajectory of $\PSD$ that
  is contained in $\MM$ and hence in $\EE(F,\xo,\yo,\zo)$.
  From~\eqref{eq:lie-v1-ineq}, we get
  \begin{align}
    & \Lie_{\PSD} V_1(x,y,z) \notag
    \\
    & \le - (x-\xo)^\top \gradient_x F(x,y,z) + (y-\yo)^\top
    \gradient_y F(x,y,z)\notag
    \\
    & \qquad + (z - \zo)^\top \gradient_z F(x,y,z) \notag
    \\
    & \le F(x,y,z) - F(x,\yo,\zo) + F(\xo,y,z) - F(x,y,z) \notag
    \\
    & \le F(\xo,\yo,\zo) - F(x,\yo,\zo) + F(\xo,y,z) \notag
    \\
    & \qquad - F(\xo,\yo,\zo) \le 0, \label{eq:V1-last-ineq}
  \end{align}
  where in the second inequality we have used the first-order
  convexity and concavity property of the maps $x \mapsto F(x,y,z)$
  and $(y,z) \mapsto F(x,y,z)$. Now since $\EE(F,\xo,\yo,\zo) =
  \setdef{(\xo,y,z)}{\Lie_{\PSD} V_1(\xo,y,z) = 0}$, using the above
  inequality, we get $F(\xo,y(t),z(t)) = F(\xo,\yo,\zo)$ for all $t
  \ge 0$.  Thus, for all $t \ge 0$, $\Lie_{\PSD} F(\xo,y(t),z(t)) = 0$
  which yields 
  \begin{multline*}
    \gradient_y F(\xo,y(t),z(t))^\top [\gradient_y
    F(\xo,y(t),z(t))]_{y(t)}^+ 
    \\ + \norm{\gradient_z F(\xo,y(t),z(t))}^2
    = 0
  \end{multline*}
  Note that both terms in the above expression are nonnegative and so,
  we get $[\gradient_y F(\xo,y(t),z(t))]_{y(t)}^+ = 0$ and
  $\gradient_z F(\xo,y(t),z(t)) = 0$ for all $t \ge 0$. In particular,
  this holds at $t = 0$ and so, $(x,y,z) \in \sdl{F}$, and we conclude
  $\MM \subset \sdl{F}$. Hence $\sdl{F}$ is globally asymptotically
  stable.  Combining this with the fact that individual saddle points
  are stable, one deduces the pointwise convergence of trajectories
  along the same lines as in~\cite[Corollary~5.2]{SPB-DSB:03}.
\end{IEEEproof}

A closer look at the proof of the above result reveals that the same
conclusion also holds under milder conditions on the saddle
function. In particular, $F$ need only be twice continuously
differentiable in a neighborhood of the saddle point and the local
strong convexity-concavity can be relaxed to a condition on the line
integral of Hessian blocks of $F$.  We state next this stronger
result. 

\begin{theorem}\longthmtitle{Global asymptotic stability of the set of
    saddle points under~$\PSD$}\label{th:proj-local2global-2}
  Let $F$ be convex-concave and continuously differentiable with
  locally Lipschitz gradient. Suppose there is a saddle point
  $(\xo,\yo,\zo)$ and a neighborhood of this point $\UU_* \subset
  \real^n \times \realnonnegative^p \times \real^m$ such that $F$ is
  twice continuously differentiable on $\UU_*$ and either of the
  following holds
  \begin{enumerate}
    \item for all $(x,y,z) \in \UU_*$, 
      \begin{align*}
        \int_0^1 \gradient_{xx} F(x(s),y(s),z(s)) ds \succ 0,
      \end{align*}
    \item for all $(x,y,z) \in \UU_*$, 
      \begin{align*}
	\int_0^1 \left[\begin{array}{c c}\gradient_{yy} F &
            \gradient_{yz} F
            \\
            \gradient_{zy} F & \gradient_{zz} F
       \end{array}\right]_{(x(s),y(s),z(s))} ds \prec 0,
     \end{align*}
  \end{enumerate}
  where $(x(s),y(s),z(s))$ are given in~\eqref{eq:Hdefs}. 
  Then, $\sdl{F}$ is globally asymptotically stable
  under the projected saddle-point dynamics $\PSD$ and the convergence
  of trajectories is to a point.
\end{theorem}

We omit the proof of this result for space reasons: the argument is
analogous to the proof of Theorem~\ref{th:proj-local2global}, where
one replaces the integral of Hessian blocks by the integral of
generalized Hessian blocks (see~\cite[Chapter 2]{FHC:83} for the
definition of the latter), as the function is not twice continuously
differentiable everywhere.

\begin{figure}[htb]
  \centering%
  \subfloat[$(x,y,z)$]{\includegraphics[width = 0.48
    \linewidth]{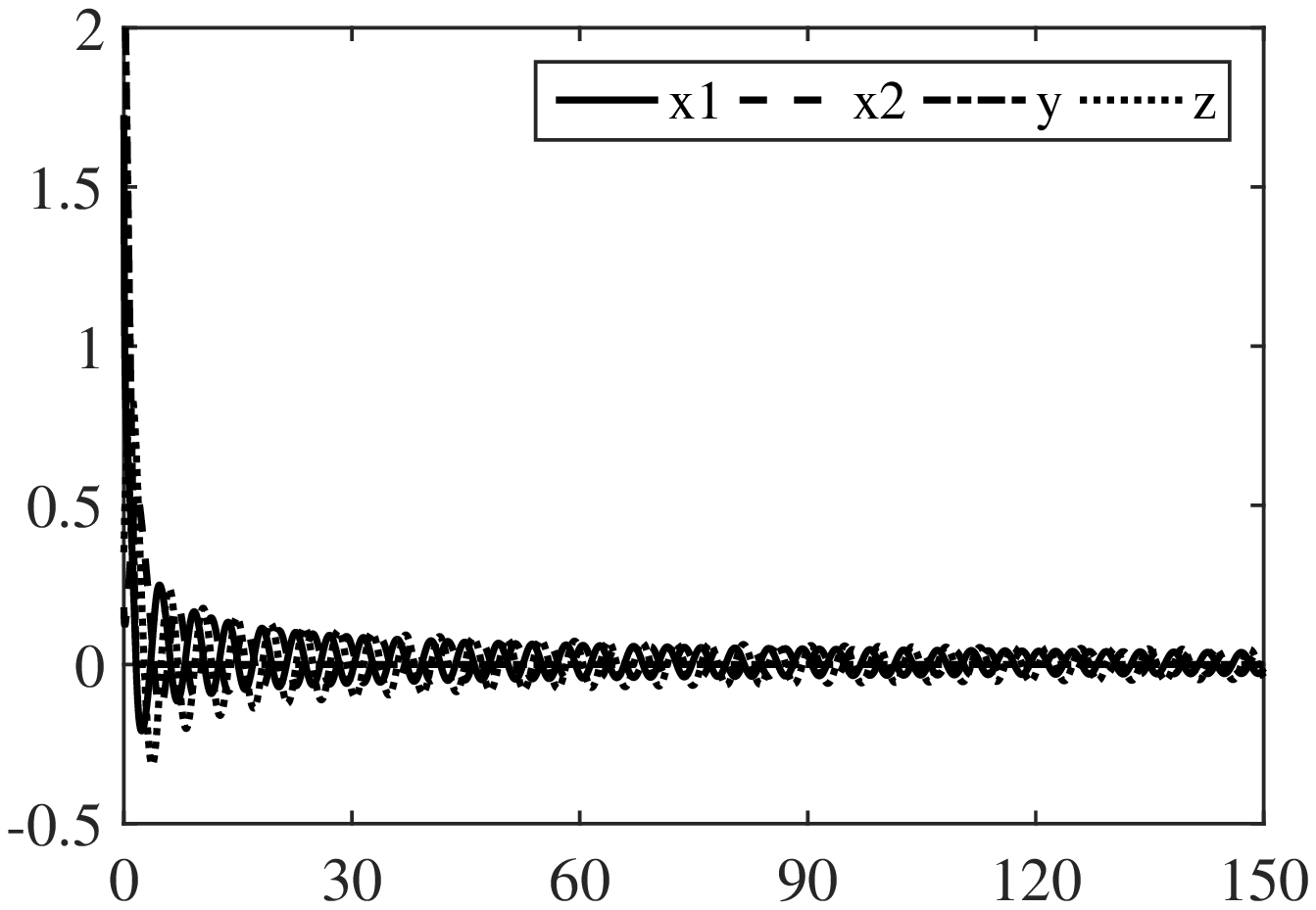}}
  \, \,
  \subfloat[$V_1$]{\includegraphics[width = 0.457
    \linewidth]{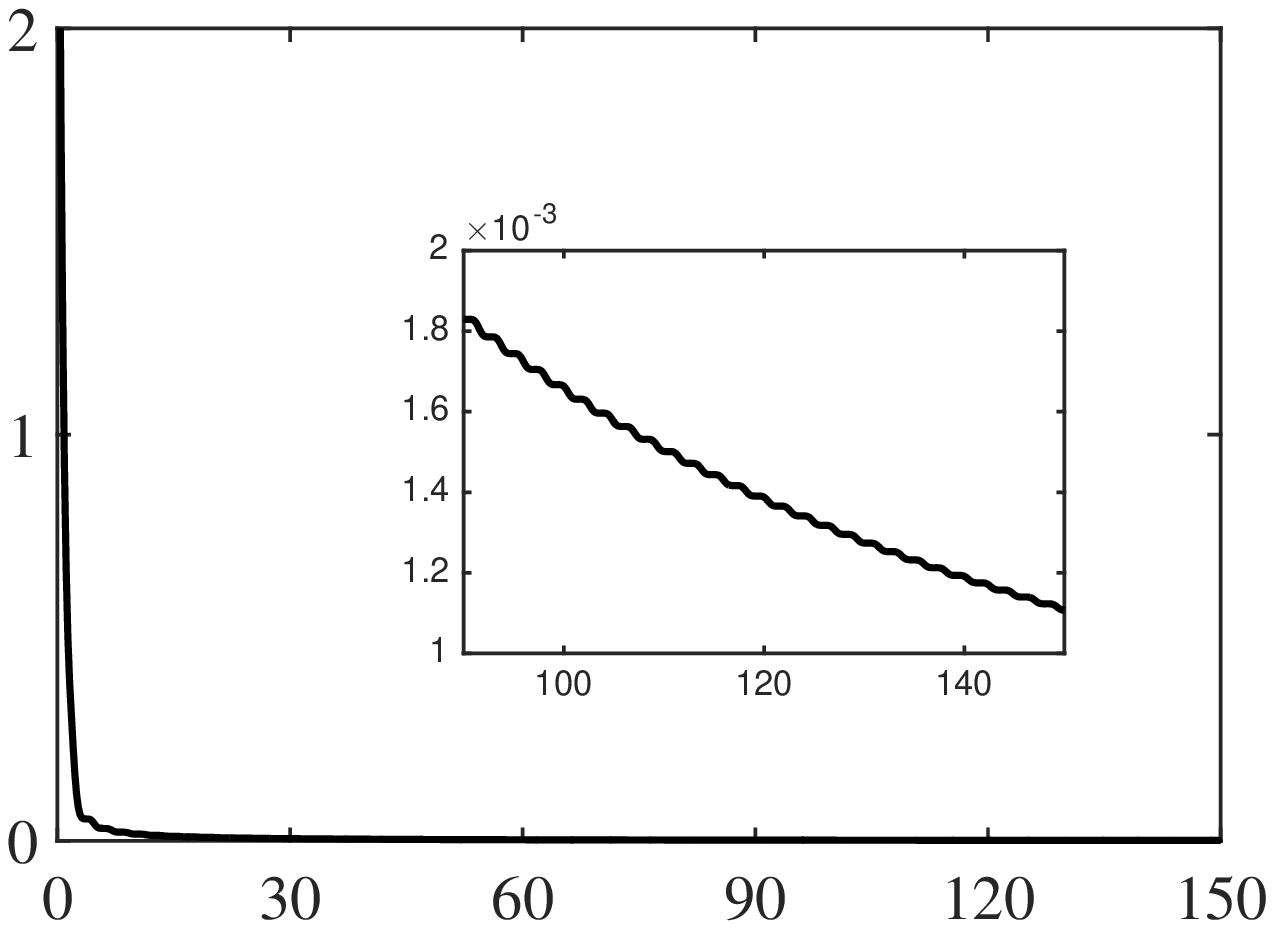}}
  \caption{Execution of the projected saddle-point
    dynamics~\eqref{eq:proj-saddle-dyn} starting from $(1.7256,
    0.1793, 2.4696, 0.3532)$ for Example~\ref{ex:one}.  As guaranteed
    by Theorem~\ref{th:proj-local2global-2}, the trajectory converges
    to the unique saddle point $0$ and the function $V_1$ defined
    in~\eqref{eq:V1} decreases
    monotonically.}\label{fig:local-to-global}
\end{figure}

\begin{example}\longthmtitle{Illustration of global asymptotic
    convergence}\label{ex:one}
  {\rm Consider $\map{F}{\real^2 \times \realnonnegative \times
      \real}{\real}$ given as
    \begin{align}\label{eq:ex-1-F}
      F(x,y,z) = f(x) +
      y(-x_1 - 1) + z(x_1 - x_2),
    \end{align}
    where
    \begin{align*}
      f(x) = \begin{cases} \norm{x}^4, & \quad \text{ if } \norm{x}
        \le \frac{1}{2},
        \\
        \frac{1}{16} + \frac{1}{2}(\norm{x} - \frac{1}{2}), & \quad
        \text{ if } \norm{x} \ge
        \frac{1}{2}.  \end{cases} \end{align*} Note that $F$ is
    convex-concave on $(\real^2) \times (\realnonnegative \times
    \real)$ and $\sdl{F} = \{0\}$. Also, $F$ is continuously
    differentiable on the entire domain and its gradient is locally
    Lipschitz. Finally, $F$ is twice continuously differentiable on
    the neighborhood $\UU_* = B_{1/2}(0) \cap (\real^2 \times
    \realnonnegative \times \real)$ of the saddle point $0$ and
    hypothesis (i) of Theorem~\ref{th:proj-local2global-2} holds on~
    $\UU_*$. Therefore, we conclude from
    Theorem~\ref{th:proj-local2global-2} that the trajectories of the
    projected saddle-point dynamics of $F$ converge globally
    asymptotically to the saddle point $0$.
    Figure~\ref{fig:local-to-global} shows an execution.  } \oprocend
\end{example}

\begin{remark}\longthmtitle{Comparison with the
    literature}\label{re:connect-lit}
  {\rm Theorems~\ref{th:proj-local2global}
    and~\ref{th:proj-local2global-2} complement the available results
    in the literature concerning the asymptotic convergence properties
    of
    saddle-point~\cite{KA-LH-HU:58,AC-BG-JC:17-sicon,SKN-JC:16-sicon}
    and primal-dual dynamics~\cite{DF-FP:10,AC-EM-JC:16-scl}. The
    former dynamics corresponds to~\eqref{eq:proj-saddle-dyn} when the
    variable $y$ is absent and the later to~\eqref{eq:proj-saddle-dyn}
    when the variable $z$ is absent. For both saddle-point and
    primal-dual dynamics, existing
    global asymptotic stability results require assumptions on the
    global properties of $F$, in addition to the global
    convexity-concavity of $F$, such as global strong
    convexity-concavity~\cite{KA-LH-HU:58}, global strict
    convexity-concavity, and its
    generalizations~\cite{AC-BG-JC:17-sicon}.  In contrast, the
    novelty of our results lies in establishing that certain local
    properties of the saddle function are enough to guarantee global
    asymptotic convergence.  \oprocend }
\end{remark}

\section{Lyapunov function for constrained convex optimization
  problems}\label{sec:Lyapunov-function}

Our discussion above has established the global asymptotic stability
of the set of saddle points resorting to LaSalle-type arguments
(because the function $V_1$ defined in~\eqref{eq:V1} is not a strict
Lyapunov function).  In this section, we identify instead a strict
Lyapunov function for the projected saddle-point dynamics when the
saddle function~$F$ corresponds to the Lagrangian of a constrained
optimization problem, \changes{cf. Remark~\ref{re:motivation}}.  The
relevance of this result stems from two facts. On the one hand, the
projected saddle-point dynamics has been employed profusely to solve
network optimization problems. On the other hand, although the
conclusions on the asymptotic convergence of this dynamics that can be
obtained with the identified Lyapunov function are the same as in the
previous section, having a Lyapunov function available is advantageous
for a number of reasons, including the study of robustness against
disturbances, the characterization of the algorithm convergence rate,
or as a design tool for developing opportunistic state-triggered
implementations. We come back to this point \changes{in
  Section~\ref{sec:ISS+self} below}.

\begin{theorem}\longthmtitle{Lyapunov function for
    $\PSD$}\label{th:lyapunov-conv-4}
  Let $\map{F}{\real^n \times \realnonnegative^p \times
    \real^m}{\real}$ be defined as
  \begin{align}\label{eq:ineq-eq-const-opt}
    F(x,y,z) = f(x) + y^\top g(x) + z^\top (Ax - b),
  \end{align}
  where $\map{f}{\real^n}{\real}$ is strongly convex, twice
  continuously differentiable, $\map{g}{\real^n}{\real^p}$ is convex,
  twice continuously differentiable, $A \in \real^{m \times n}$, and
  $b \in \real^m$. For each $(x,y,z) \in \real^n \times
  \realnonnegative^p \times \real^m$, define the index set of active
  constraints
  \begin{multline*}
    \JJ(x,y,z)  = \setdef{j \in \until{p}}{y_j = 0 \text{ and }
     \\
     (\gradient_y F(x,y,z))_j < 0}.
  \end{multline*}
  Then, the function $\map{V_2}{\real^n \times \realnonnegative^p
  \times \real^m}{\real}$,
  \begin{multline*}
    V_2(x,y,z) = \frac{1}{2} \Bigl(\norm{\gradient_x F(x,y,z)}^2 +
    \norm{\gradient_z F(x,y,z)}^2
    \\
    + \sum_{j \in \until{p} \setminus \JJ(x,y,z)} ((\gradient_y
    F(x,y,z))_j)^2 \Bigr)
    \\
    + \frac{1}{2}\norm{(x,y,z)}_{\sdl{F}}^2
  \end{multline*}
  is nonnegative everywhere in its domain and $V_2(x,y,z) = 0$ if and
  only if $(x,y,z) \in \sdl{F}$. \changes{Moreover, for any trajectory
    $t \mapsto (x(t),y(t),z(t))$ of $\PSD$, the map $t \mapsto
    V_2(x(t),y(t),z(t))$ 
  \begin{enumerate}
  \item is differentiable almost everywhere and if $(x(t),y(t),z(t))
    \not \in \sdl{F}$ for some $t \ge 0$, then $\frac{d}{dt}
    V_2(x(t),y(t),z(t)) < 0$ provided the derivative exists.
    Furthermore, for any sequence of times $\{t_k\}_{k=1}^\infty$ such
    that $t_k \to t$ and $\frac{d}{dt} V_2(x(t_k),y(t_k),z(t_k))$
    exists for every $t_k$, we have $\limsup_{k \to \infty}
    \frac{d}{dt} V(x(t_k),y(t_k),z(t_k)) < 0$,
  \item is right-continuous and at any point of discontinuity $t' \ge
    0$, we have $V_2(x(t'),y(t'),z(t')) \le \lim_{t \uparrow t'}
    V_2(x(t),y(t),z(t))$.
  \end{enumerate}
  }
  As a consequence, $\sdl{F}$ is globally asymptotically
  stable under $\PSD$ and convergence of trajectories is to a point.
\end{theorem}
\begin{IEEEproof}
  We start by partitioning the domain based on the active constraints.
  Let $\II \subset \until{p}$ and
  \begin{align*}
    \DD(\II) = \setdef{(x,y,z) \in \real^n \times \realnonnegative^p
      \times \real^m}{\JJ(x,y,z) = \II}.
  \end{align*}
  Note that for $\II_1, \II_2 \subset \until{p}$, $\II_1 \not =
  \II_2$, we have $\DD(\II_1) \cap \DD(\II_2) = \emptyset$. Moreover,  
  \begin{align*}
    \real^n \times \realnonnegative^p \times \real^m = \bigcup_{\II
      \subset \until{p}} \DD(\II).
  \end{align*}
  For each $\II \subset \until{p}$, define the function
  \begin{multline}\label{eq:lyap-set}
    V_2^{\II}(x,y,z) = \frac{1}{2} \Bigl(\norm{\gradient_x
      F(x,y,z)}^2 + \norm{\gradient_z F(x,y,z)}^2
    \\
    + \sum_{j \not \in \II} ((\gradient_y F(x,y,z))_j)^2 \Bigr) +
    \frac{1}{2}\norm{(x,y,z)}_{\sdl{F}}^2.
  \end{multline}
  \changes{These functions will be used later for analyzing the
    evolution of $V_2$.  Consider a trajectory $t \mapsto
    (x(t),y(t),z(t))$ of $\PSD$ starting at some point
    $(x(0),y(0),z(0)) \in \real^n \times \realnonnegative^p \times
    \real^m$. Our proof strategy consists of proving assertions (i)
    and (ii) for two scenarios, depending on whether or not there
    exists $\delta > 0$ such that the difference between two
    consecutive time instants when the trajectory switches from one
    partition set to another is lower bounded by~$\delta$.
    
    \textbf{Scenario 1: time elapsed between consecutive switches is
      lower bounded:}
    Let $(a,b) \subset \realnonnegative$, $b - a \ge \delta$, be a
    time interval for which the trajectory belongs to a partition
    $\DD(\II')$, $\II' \subset \until{p}$, for all $t \in (a,b)$. In
    the following, we show that $\frac{d}{dt} V_2(x(t),y(t),z(t))$
    exists for almost all $t \in (a,b)$ and its value is negative
    whenever $(x(t),y(t),z(t)) \not \in \sdl{F}$. Consider the
    function $V_2^{\II'}$ defined in~\eqref{eq:lyap-set}} and note that
  $t \mapsto V_2^{\II'}(x(t),y(t),z(t))$ is absolutely continuous as
  $V_2^{\II'}$ is continuously differentiable on $\real^n \times
  \realnonnegative^p \times \real^m$ and the trajectory is
  absolutely continuous.
  Employing Rademacher's Theorem~\cite{FHC:83}, we deduce that the map
  $t \mapsto V_2^{\II'}(x(t),y(t),z(t))$ is differentiable almost
  everywhere.  By definition, $V_2(x(t),y(t),z(t)) =
  V_2^{\II'}(x(t),y(t),z(t))$ for all $t \in (a,b)$.  Therefore
  \begin{align}\label{eq:first-equality}
    \frac{d}{dt} V_2(x(t),y(t),z(t)) = \frac{d}{dt}
    V_2^{\II'}(x(t),y(t),z(t))
  \end{align}
  for almost all $t \in (a,b)$. Further, since $V_2^{\II'}$ is continuously
  differentiable, we have
  \begin{align}\label{eq:second-equality}
    \frac{d}{dt} V_2^{\II'}(x(t),y(t),z(t)) = \Lie_{\PSD}
    V_2^{\II'}(x(t),y(t),z(t)) .
  \end{align}
  Now consider any $(x,y,z) \in \DD(\II') \setminus \sdl{F}$. Our next
  computation shows that $\Lie_{\PSD}V_2^{\II'}(x,y,z) < 0$. We have
  \begin{align}
    & \Lie_{\PSD}V_2^{\II'}(x,y,z) \notag
    \\
    & = - \gradient_x F(x,y,z)^\top \gradient_{xx} F(x,y,z) \gradient_x
    F(x,y,z) \notag
    \\
    & \quad + \begin{bmatrix} [\gradient_y F(x,y,z)]_y^+
      \\
      \gradient_z F(x,y,z)
    \end{bmatrix}^\top
    \begin{bmatrix}
      \gradient_{yy} F & \gradient_{yz} F
      \\
      \gradient_{zy} F & \gradient_{zz} F \end{bmatrix}_{(x,y,z)} \notag
    \\
    & \qquad \qquad \qquad \qquad \qquad \qquad \qquad \qquad
    \begin{bmatrix}
      [\gradient_y F(x,y,z)]_y^+
      \\
      \gradient_z F(x,y,z)
    \end{bmatrix} \notag
    \\
    & \quad + \Lie_{\PSD} \Bigl(\frac{1}{2} \norm{(x,y,z)}_{\sdl{F}}^2
    \Bigr). \label{eq:vI-ineq}
  \end{align} 
  The first two terms in the above expression are the Lie derivative
  of $(x,y,z) \mapsto V_2^{\II'}(x,y,z) - \frac{1}{2}
  \norm{(x,y,z)}_{\sdl{F}}^2$. This computation can be shown using the
  properties of the operator $[\cdot]_y^+$. Now let $(\xo,\yo,\zo) =
  \proj_{\sdl{F}}(x,y,z)$. Then, by Danskin's Theorem~\cite[p.
  99]{FHC-YSL-RJS-PRW:98}, we have
  \begin{equation}\label{eq:danskin-1}
    \gradient \norm{(x,y,z)}_{\sdl{F}}^2 = 2(x-\xo;y-\yo;z-\zo)
  \end{equation}
  Using this expression, we get
  \begin{align*}
    & \Lie_{\PSD} \Bigl(\frac{1}{2} \norm{(x,y,z)}_{\sdl{F}}^2
    \Bigr) 
    \\
    & = -(x-\xo)^\top \gradient_x F(x,y,z) + (y-\yo)^\top
    [\gradient_y F(x,y,z)]_y^+ 
    \\
    & \quad + (z-\zo)^\top \gradient_z F(x,y,z)
    \\
    & \le F(\xo,y,z) - F(\xo,\yo,\zo) + F(\xo,\yo,\zo)
    \\
     & \quad - F(x,\yo,\zo),
   \end{align*}
  where the last inequality follows from~\eqref{eq:V1-last-ineq}. Now
  using the above expression in~\eqref{eq:vI-ineq} we get
  \begin{align*}
    & \Lie_{\PSD}V_2^{\II'}(x,y,z) 
    \\
    & \le - \gradient_x F(x,y,z) \gradient_{xx} F(x,y,z) \gradient_x
    F(x,y,z) 
    \\
    & \quad + \begin{bmatrix} [\gradient_y F(x,y,z)]_y^+
      \\
      \gradient_z F(x,y,z)
    \end{bmatrix}^\top
    \begin{bmatrix}
      \gradient_{yy} F & \gradient_{yz} F
      \\
      \gradient_{zy} F & \gradient_{zz} F \end{bmatrix}_{(x,y,z)} 
    \\
    & \qquad \qquad \qquad \qquad \qquad \qquad \qquad \qquad
    \begin{bmatrix}
      [\gradient_y F(x,y,z)]_y^+
      \\
      \gradient_z F(x,y,z)
    \end{bmatrix} 
    \\
    & \quad  +  F(\xo,y,z) - F(\xo,\yo,\zo) + F(\xo,\yo,\zo)
    \\
    & \quad - F(x,\yo,\zo) \le 0.
  \end{align*} 
  If $\Lie_{\PSD} V_2^{\II'}(x,y,z) = 0$, then (a) $\gradient_x
  F(x,y,z) = 0$; (b) $x = \xo$; and (c) $F(\xo,y,z) =
  F(\xo,\yo,\zo)$. From (b) and~\eqref{eq:equi-cond}, we conclude that
  $\gradient_z F(x,y,z) = 0$.  From (c)
  and~\eqref{eq:ineq-eq-const-opt}, we deduce that $(y - \yo)^\top
  g(\xo) =0$. Note that for each $i \in \until{p}$, we have $(y_i -
  (\yo)_i) (g(\xo))_i \le 0$. This is because either $(g(\xo))_i = 0$
  in which case it is trivial or $(g(\xo))_i < 0$ in which case
  $(\yo)_i = 0$ (as $\yo$ maximizes the map $y \mapsto y^\top g(\xo)$)
  thereby making $y_i - (\yo)_i \ge 0$. Since, $(y_i - (\yo)_i)
  (g(\xo))_i \le 0$ for each $i$ and $(y - \yo)^\top g(\xo) =0$, we
  get that for each $i \in \until{p}$, either $(g(\xo))_i = 0$ or $y_i
  = (\yo)_i$. Thus, $[\gradient_y F(x,y,z)]_y^+ = 0$. These facts
  imply that $(x,y,z) \in \sdl{F}$. Therefore, if $(x,y,z) \in
  \DD(\II')\setminus \sdl{F}$ then $\Lie_{\PSD} V_2^{\II'}(x,y,z)< 0$.
  \changes{Combining this with~\eqref{eq:first-equality}
    and~\eqref{eq:second-equality}, we deduce
    \begin{align*}
      \frac{d}{dt} V_2(x(t),y(t),z(t)) < 0
    \end{align*}
    for almost all $t \in (a,b)$.  Therefore, between any two switches
    in the partition, the evolution of $V_2$ is differentiable and the
    value of the derivative is negative. Since the number of time
    instances when a switch occurs is countable, the first part of
    assertion (i) holds. To show the limit condition, consider $t \ge
    0$ such that $(x(t),y(t),z(t)) \not \in \sdl{F}$. Let
    $\{t_k\}_{k=1}^\infty$ be such that $t_k \to t$ and $\frac{d}{dt}
    V_2(x(t_k),y(t_k),z(t_k))$ exists for every $t_k$. By continuity,
    $\lim_{k \to \infty} (x(t_k),y(t_k),z(t_k)) =
    (x(t),y(t),z(t))$. Let $\BB \subset \real^n \times
    \realnonnegative^p \times \real^m$ be a compact neighborhood of
    $(x(t),y(t),z(t))$ such that $\BB \cap \sdl{F} =
    \emptyset$. Without loss of generality, assume that
    $\{x(t_k),y(t_k),z(t_k))\}_{k=1}^\infty \subset \BB$. Define
    \begin{align*}
      S = \max \setdef{\Lie_{\PSD} V_2^{\JJ(x,y,z)} (x,y,z)}{(x,y,z)
        \in \BB}.
    \end{align*}
    The Lie derivatives in the above expression are well-defined and
    continuous as each $V_2^{\JJ(x,y,z)}$ is continuously
    differentiable.  Note that $S < 0$ as $\BB \cap \sdl{F} =
    \emptyset$. Moreover, as established above, for each $k$,
    $\frac{d}{dt} V_2(x(t_k),y(t_k),z(t_k)) = \Lie_{\PSD}
    V_2^{\JJ(x(t_k),y(t_k),z(t_k))} (x(t_k),y(t_k),z(t_k)) \le S$.
    Thus, we get $\limsup_{k \to \infty} \frac{d}{dt}
    V_2(x(t_k),y(t_k),z(t_k)) \le S < 0$, establishing (i) for
    Scenario~1.}
  
  To prove assertion (ii), note that discontinuity in $V_2$ can only
  happen when the trajectory switches the partition. In order to
  analyze this, consider any time instant $t' \ge 0$ and let
  $(x(t'),y(t'),z(t')) \in \DD(\II')$ for some $\II' \subset
  \until{p}$.  \changes{Looking at times $t \ge t'$, two cases arise:
    \begin{enumerate}[label=(\alph*)]
    \item There exists $\tilde{\delta} > 0$ such that
      $(x(t),y(t),z(t)) \in \DD(\II')$ for all $t \in [t',t'+
      \tilde{\delta})$.
    \item There exists $\tilde{\delta} > 0$ and $\II \not = \II'$ such
      that $(x(t),y(t),z(t)) \in \DD(\II)$ for all $t \in
      (t',t'+\tilde{\delta})$.
    \end{enumerate}
    One can show that for Scenario 1, the trajectory cannot show any
    behavior other than the above mentioned two cases.  We proceed to
    show that in both the above outlined cases, $t \mapsto
    V_2(x(t),y(t),z(t))$ is right-continuous at $t'$.  Case (a) is
    straightforward as $V_2$ is continuous in the domain $\DD(\II')$
    and the trajectory is absolutely continuous.  In case (b), $\II
    \not = \II'$ implies that there exists $j \in \until{p}$ such that
    either $j \in \II \setminus \II'$ or $j \in \II' \setminus
    \II$. Note that the later scenario, i.e., $j \in \II'$ and $j \not
    \in \II$ cannot happen.  Indeed by definition $(y(t'))_j = 0$ and
    $(\gradient_y F(x(t'),y(t'),z(t')))_j < 0$ and by continuity of
    the trajectory and the map $\gradient_y F$, these conditions also
    hold for some finite time interval starting at $t'$.  Therefore,
    we focus on the case that $j \in \II \setminus \II'$.  Then,
    either $(y(t'))_j > 0$ or $(\gradient_y F(x(t'),y(t'),z(t')))_j
    \ge 0$. The former implies, due to continuity of trajectories,
    that it is not possible to have $j \in \II$. Similarly, by
    continuity if $(\gradient_y F(x(t'),y(t'),z(t')))_j > 0$, then one
    cannot have $j \in \II$.  Therefore, the only possibility is
    $(y(t'))_j = 0$ and $(\gradient_y F(x(t'),y(t'),z(t')))_j = 0$.
    This implies that the term $t \mapsto (\gradient_y
    F(x(t),y(t),z(t)))^2_j$ is right-continuous at $t'$.  Since this
    holds for any $j \in \II \setminus \II'$, we conclude
    right-continuity of $V_2$ at $t'$. Therefore, for both cases (a)
    and (b), we conclude right-continuity of $V_2$.
  
    Next we show the limit condition of assertion (ii). Let $t' \ge 0$
    be a point of discontinuity.} Then, from the preceding discussion,
  there must exist $\II, \II' \subset \until{p}$, $\II \not = \II'$,
  such that $(x(t'),y(t'),z(t')) \in \DD(\II')$ and $(x(t),y(t),z(t))
  \in \DD(\II)$ for all $t \in (t'-\delta, t')$. By continuity,
  $\lim_{t \uparrow t'} V_2(x(t),y(t),z(t))$ exists. Note that if $j
  \in \II$ and $j \not \in \II'$, then the term getting added to $V_2$
  at time $t'$ which was absent at times $t \in (t'-\delta,t')$, i.e.,
  $(\gradient_y F(x(t),y(t),z(t)))_j^2$, is zero at $t'$. Therefore,
  the discontinuity at $t'$ can only happen due to the existence of $j
  \in \II' \setminus \II$.  That is, a constraint becomes active at
  time $t'$ which was inactive in the time interval $(t'-\delta,t')$.
  Thus, the function $V_2$ loses a nonnegative term at time $t'$. This
  can only mean at $t'$ the value of $V_2$ decreases. Hence, the limit
  condition of assertion (ii) holds.

  \changes{\textbf{Scenario 2: time elapsed between consecutive
      switches is not lower bounded:} Observe that three cases arise.
    First is when there are only a finite number of switches in
    partition in any compact time interval.
    In this case, the analysis of Secnario 1 applies to every compact
    time interval and so assertions (i) and (ii) hold.  The second
    case is when there exist time instants $t' > 0$ where there is
    absence of ``finite dwell time'', that is, there exist index sets
    $\II_1 \not = \II_2$ and $\II_2 \not = \II_3$ such that
    $(x(t),y(t),z(t)) \in \DD(\II_1)$ for all $t \in (t'-\eps_1,t')$
    and some $\eps_1 > 0$; $(x(t'),y(t'),z(t')) \in \DD(\II_2)$; and
    $(x(t),y(t),z(t)) \in \DD(\II_3)$ for all $t \in (t',t'+\eps_2)$
    and some $\eps_2 > 0$. Again using the arguments of Scenario 1,
    one can show that both assertions (i) and (ii) hold for this case
    if there is no accumulation point of such time instants~$t'$.

    The third case instead is when there are infinite switches in a
    finite time interval. We analyze this case in parts. Assume that
    there exists a sequence of times $\{t_k\}_{k=1}^\infty$, $t_k
    \uparrow t'$, such that trajectory switches partition at
    each~$t_k$. The aim is to show left-continuity of $t \mapsto
    V(x(t),y(t),z(t))$ at $t'$.
    Let $\II^s \subset \until{p}$ be the set of indices that switch
    between being active and inactive an infinite number of times
    along the sequence $\{t_k\}$ (note that the set is nonempty as
    there are an infinite number of switches and a finite number of
    indices).  To analyze the left-continuity at $t'$, we only need to
    study the possible occurrence of discontinuity due to terms in
    $V_2$ corresponding to the indices in $\II^s$, since all other
    terms do not affect the continuity.  Pick any $j \in \II^s$. Then,
    the term in~$V_2$ corresponding to the index $j$ satisfies
    \begin{align}\label{eq:switch-limit}
      \lim_{k \to \infty} (\gradient_y F(x(t_k),y(t_k),z(t_k)))_j^2 =
      0.
    \end{align}
    In order to show this, assume the contrary. This implies the
    existence of $\eps > 0$ such that
    \begin{align*}
      \liminf_{k \to \infty} (\gradient_y F(x(t_k),y(t_k),z(t_k)))_j^2
      \ge \eps.
    \end{align*}
    As a consequence, the set of $k$ for which $(\gradient_y
    F(x(t_k),y(t_k),z(t_k)))_j^2 \ge \eps/2$ is infinite. Recall that
    if the constraint $j$ becomes active at $t_k$, then $V_2$
    decreases by at least $(\gradient_y F(x(t_k),y(t_k),z(t_k)))_j^2$
    at $t_k$. Further, $V_2$ decreases montonically between any
    consecutive $t_k$'s. These facts lead to the conclusion that $V_2$
    tends to $-\infty$ as $t_k \to t'$. However, $V_2$ takes
    nonnegative values, yielding a contradiction.
    Hence,~\eqref{eq:switch-limit} is true for all $j \in \II^s$ and
    so,
    \begin{align*}
      \lim_{k \to \infty} V_2(x(t_k),y(t_k),z(t_k)) =
      V_2(x(t'),y(t'),z(t')),
    \end{align*}
    proving left-continuity of $V_2$ at $t'$.  Using this reasoning,
    one can also conclude that if the infinite number of switches
    happen on a sequence $\{t_k\}_{k=1}^\infty$ with $t_k \downarrow
    t'$, then one has right-continuity at $t'$.  Therefore, at each
    time instant when a switch happens, we have right-continuity of $t
    \mapsto V_2(x(t),y(t),z(t))$ and at points where there is
    accumulation of switches we have continuity (depending on which
    side of the time instance the accumulation takes place).  This
    proves assertion (ii).  Note that in this case too we have a
    countable number of time instants where the partition set switches
    and so the map $t \mapsto V_2(x(t),y(t),z(t))$ is differentiable
    almost everywhere. Moreover, one can also analyze, as done in
    Scenario 1, that the limit condition of assertion (i) holds in
    this case. These facts together establish the condition of
    assertion (ii), completing the proof.  }
\end{IEEEproof}

\begin{remark}\longthmtitle{Multiple  Lyapunov functions}\label{re:charac-V2}
  {\rm The Lyapunov function $V_2$ is discontinuous on the domain
    $\real^n \times \realnonnegative^p \times \real^m$. However, it
    can be seen as multiple (continuously differentiable) Lyapunov
    functions~\cite{MSB:98}, each valid on a domain, patched together in
    an appropriate way such that along the trajectories of $\PSD$, the
    evolution of $V_2$ is continuously differentiable with negative
    derivative at intervals where it is continuous and at times of
    discontinuity the value of $V_2$ only decreases. Note that in the
    absence of the projection in $\PSD$ (that is, no $y$-component of
    the dynamics), the function $V_2$ takes a much simpler
    form with no discontinuities and is continuously differentiable on the
    entire domain.
  }
  \oprocend
\end{remark}

\begin{remark}\longthmtitle{Connection with the literature: II}
  {\rm The two functions whose sum defines $V_2$ are, individually by
    themselves, sufficient to establish asymptotic convergence
    of~$\PSD$ using LaSalle Invariance arguments, see
    e.g.,~\cite{DF-FP:10,AC-EM-JC:16-scl}. However, the fact that
    their combination results in a strict Lyapunov function for the
    projected saddle-point dynamics is a novelty of our analysis here.
    In~\cite{SKN-JC:16-sicon}, a different Lyapunov function is
    proposed and an exponential rate of convergence is established for
    a saddle-point-like dynamics which is similar to $\PSD$ but
    without projection components.  } \oprocend
\end{remark}

\section{ISS and self-triggered implementation of the saddle-point
  dynamics}\label{sec:ISS+self}

Here, we build on the novel Lyapunov function identified in
Section~\ref{sec:Lyapunov-function} to explore other properties of the
projected saddle-point dynamics beyond global asymptotic convergence.
Throughout this section, we consider saddle functions~$F$ that
corresponds to the Lagrangian of an equality-constrained optimization
problem, i.e.,
\begin{align}\label{eq:F-wo-y}
  F(x,z) = f(x) + z^\top (Ax - b),
\end{align}
where $A \in \real^{m \times n}$, $b \in \real^m$, and
$\map{f}{\real^n}{\real}$. The reason behind this focus is that, in
this case, the dynamics~\eqref{eq:proj-saddle-dyn} is smooth and the
Lyapunov function identified in Theorem~\ref{th:lyapunov-conv-4} is
continuously differentiable.  These simplifications allow us to
analyze input-to-state stability of the dynamics using the theory of
ISS-Lyapunov functions (cf. Section~\ref{subsec:iss-prelim}).  On the
other hand, we do not know of such a theory for projected systems,
which precludes us from carrying out ISS analysis for
dynamics~\eqref{eq:proj-saddle-dyn} for a general saddle function.
The projected saddle-point dynamics~\eqref{eq:proj-saddle-dyn} for the
class of saddle functions given in~\eqref{eq:F-wo-y} takes the form
\begin{subequations}\label{eq:saddle-dyn}
  \begin{align}
    \dot x & = - \gradient_x F(x,z) = - \gradient f(x) - A^\top z,
    \\
    \dot z & = \gradient_z F(x,z) = Ax - b,
  \end{align}
\end{subequations}
corresponding to equations~\eqref{eq:proj-saddle-dyn-1}
and~\eqref{eq:proj-saddle-dyn-3}. We term these dynamics simply
\emph{saddle-point dynamics} and denote it as $\map{\SD}{\real^n
  \times \real^m}{\real^n \times \real^m}$.

\subsection{Input-to-state stability}\label{subsec:iss}

Here, we establish that the saddle-point
dynamics~\eqref{eq:saddle-dyn} is ISS with respect to the set
$\sdl{F}$ when disturbance inputs affect it additively.  Disturbance
inputs can arise when implementing the saddle-point dynamics as a
controller of a physical system because of a variety of malfunctions,
including errors in the gradient computation, noise in state
measurements, and errors in the controller implementation. In such
scenarios, the following result shows that the
dynamics~\eqref{eq:saddle-dyn} exhibits a graceful degradation of its
convergence properties, one that scales with the size of the
disturbance.

\begin{theorem}\longthmtitle{ISS of saddle-point
    dynamics}\label{th:iss-lyapunov}
  Let the saddle function $F$ be of the form~\eqref{eq:F-wo-y}, with $f$
  strongly convex, twice continuously differentiable, and satisfying $mI
  \preceq \gradient^2 f(x) \preceq MI$ for all $x \in \real^n$ and some
  constants $0<m \le M < \infty$.  Then, the dynamics
  \begin{align}\label{eq:noise-saddle-dyn}
    \begin{bmatrix}
      \dot x
      \\
      \dot z
    \end{bmatrix}
    =
    \begin{bmatrix}
      - \gradient_x F(x,z)
      \\ 
      \gradient_z F(x,z)
    \end{bmatrix}
    + 
    \begin{bmatrix}
      u_x
      \\
      u_z
    \end{bmatrix},
  \end{align}
  where $\map{(u_x,u_z)}{\realnonnegative}{\real^n \times \real^m}$ is
  a measurable and locally essentially bounded map, is ISS with
  respect to $\sdl{F}$.
\end{theorem}
\begin{IEEEproof}
  For notational convenience, we refer to~\eqref{eq:noise-saddle-dyn}
  by $\map{\SDp}{\real^{n} \times \real^{m} \times \real^n \times
    \real^m}{\real^n \times \real^m}$.
  Our proof consists of establishing that the function
  $\map{V_3}{\real^n \times \real^m}{\realnonnegative}$,
  \begin{multline}\label{eq:V2}
    V_3(x,z) = \frac{\beta_1}{2}\norm{\SD(x,z)}^2 + 
    \frac{\beta_2}{2}\norm{(x,z)}_{\sdl{F}}^2
  \end{multline}
  with $\beta_1 > 0$, $\beta_2 = \frac{4 \beta_1 M^4}{m^2}$, is an
  ISS-Lyapunov function with respect to $\sdl{F}$ for $\SDp$. The
  statement then directly follows from
  Proposition~\ref{pr:iss-lyap-iss}.  

  We first show~\eqref{eq:iss-lyap-c1} for $V_3$, that is, there exist
  $\alpha_1, \alpha_2 > 0$ such that $\alpha_1 \norm{(x,z)}_{\sdl{F}}^2
  \le V_3(x,z) \le \alpha_2 \norm{(x,z)}_{\sdl{F}}^2$ for all $(x,z) \in
  \real^n \times \real^m$.  The lower bound follows by choosing
  $\alpha_1 = \beta_2/2$. For the upper bound, define the function
  $\map{U}{\real^n \times \real^n}{\real^{n \times n}}$ by 
  \begin{align}\label{eq:U-function}
    U(x_1,x_2) = \int_0^1 \gradient^2 f(x_1 + s(x_2 - x_1)) ds.
  \end{align}
  By assumption, it holds that $mI \preceq U(x_1,x_2) \preceq MI$ for
  all $x_1, x_2 \in \real^n$. Also, from the fundamental theorem of
  calculus, we have $\gradient f(x_2) - \gradient f(x_1) = U(x_1,x_2)
  (x_2 - x_1)$ for all $x_1, x_2 \in \real^n$.  Now pick any $(x,z) \in
  \real^n \times \real^m$. Let $(\xo,\zo) = \proj_{\sdl{F}}(x,z)$, that
  is, the projection of $(x,z)$ on the set $\sdl{F}$. This projection is
  unique as $\sdl{F}$ is convex. Then, one can write 
  \begin{align}\label{eq:grad-exp-1}
    \gradient_x F(x,z) & = \gradient_x F(\xo,\zo) + \int_0^1
    \gradient_{xx} F(x(s),z(s)) (x - \xo) ds \notag 
    \\ 
    & \quad + \int_0^1 \gradient_{zx}
    F(x(s),z(s)) (z - \zo) ds, \notag
    \\
     & = U(\xo,x) (x - \xo) + A^\top (z - \zo),
  \end{align}
  where $x(s) = \xo + s(x-\xo)$ and $z(s) = \zo + s(z - \zo)$.
  Also, note that 
  \begin{align}\label{eq:grad-exp-2}
    \gradient_z F(x,z) & = \gradient_z F(\xo,\zo) + \int_0^1
    \gradient_{xz} F(x(s),z(s)) (x - \xo) ds \notag
    \\
    & = A(x - \xo).
  \end{align}
  The expressions~\eqref{eq:grad-exp-1} and~\eqref{eq:grad-exp-2} use
  $\gradient_x F(\xo,\zo) = 0$, $\gradient_z F(\xo,\zo) =
  0$, and $\gradient_{zx} F(x,z) = \gradient_{xz} F(x,z)^\top = A^\top$
  for all $(x,z)$. From~\eqref{eq:grad-exp-1} and~\eqref{eq:grad-exp-2},
  we get
  \begin{align*}
    \norm{\SD(x,z)}^2 & \le \tilde{\alpha}_2 (\norm{x - \xo}^2 + \norm{z -
    \zo}^2)
    \\
    & = \tilde{\alpha}_2 \norm{(x,z)}_{\sdl{F}}^2,
  \end{align*}
  where $\tilde{\alpha}_2 = \frac{3}{2}(M^2 + \norm{A}^2)$. In the above
  computation, we have used the inequality $(a+b)^2 \le 3(a^2+b^2)$ for
  any $a,b \in \real$. The above
  inequality gives the upper bound $V_3(x,z) \le \alpha_2
  \norm{(x,z)}_{\sdl{F}}^2$, where $\alpha_2 =  \frac{3 \beta_1}{2}(M^2
  + \norm{A}^2) + \frac{\beta_2}{2}$.  
  
  The next step is to show that the
  Lie derivative of $V_3$ along the dynamics $\SDp$ satisfies the ISS
  property~\eqref{eq:iss-lyap-c2}.  Again, pick any $(x,z) \in \real^n
  \times \real^m$ and let $(\xo,\zo) = \proj_{\sdl{F}}(x,z)$. Then, by
  Danskin's Theorem~\cite[p.  99]{FHC-YSL-RJS-PRW:98}, we get 
  \begin{align*}
    \gradient \norm{(x,z)}_{\sdl{F}}^2 = 2(x-\xo;z-\zo).
  \end{align*}
  Using the above expression, one can compute the Lie derivative of
  $V_3$ along the dynamics $\SDp$ as 
  \begin{align*}
     \Lie_{\SDp}& V_3(x,z)  = - \beta_1 \gradient_x
    F(x,z) \gradient_{xx} F(x,z) \gradient_x F(x,z)
    \\ 
    & \quad -\beta_2(x - \xo)^\top \gradient_x F(x,z) + \beta_2(z-\zo)^\top
    \gradient_z F(x,z)
    \\ 
    & \quad + \beta_1\gradient_x F(x,z)^\top \gradient_{xx} F(x,z) u_x 
    \\ 
    & \quad +
    \beta_1 \gradient_x F(x,z)^\top \gradient_{xz} F(x,z) u_z 
    \\ 
    & \quad + \beta_1 \gradient_z F(x,z)^\top \gradient_{zx} F(x,z) u_x 
    \\ 
    & \quad + \beta_2 (x - \xo)^\top u_x + \beta_2 (z - \zo)^\top u_z. 
  \end{align*}
  Due to the particular form of $F$, we have
  \begin{alignat*}{2}
    \gradient_x F(x,z) & = \gradient f(x) + A^\top z, & \quad
    \gradient_z F(x,z)   &= A x - b,
    \\
    \gradient_{xx} F(x,z)  & = \gradient^2 f(x),  & \quad
    \gradient_{xz} F(x,z)   & = A^\top,
    \\ 
    \gradient_{zx} F(x,z) & = A, & \quad 
    \gradient_{zz} F(x,z) & = 0.
  \end{alignat*}
  Also, $\gradient_x F(\xo,\zo) = \gradient_x f(\xo)+ A^\top \zo = 0$
  and $\gradient_z F(\xo,\zo) = A \xo - b = 0$. Substituting these
  values in the expression of $\Lie_{\SDp}
  V_3$, replacing $\gradient_x F(x,z) = \gradient_x F(x,z) -
  \gradient_x F(\xo,\zo) = \gradient f(x) - \gradient f(\xo) + A^\top
  (z - \zo) = U(\xo,x) (x - \xo) + A^\top (z - \zo)$, and simplifying,
  \begin{align*}
    & \Lie_{\SDp} V_3(x,z) = 
    \\ 
    &  - \beta_1  (U(\xo,x) (x - \xo))^\top \gradient^2 f(x)
    (U(\xo,x) (x - \xo)) 
    \\
    & - \beta_1 (z - \zo)^\top A \gradient^2 f(x) A^\top (z -
    \zo) 
    \\
    & - \beta_1 (U(\xo,x) (x - \xo))^\top \gradient^2 f(x) A^\top
    (z - \zo) 
    \\
    & - \beta_1 (z -\zo)^\top A \gradient^2 f(x) (U(\xo,x) (x -
    \xo))
    \\
    & - (x -\xo)^\top U(\xo,x) (x - \xo)
    \\ 
    & + \beta_1 (U(\xo,x) (x - \xo) + A^\top (z - \zo))^\top
    \gradient^2 f(x)  u_x
    \\ 
    & + \beta_1 (U(\xo,x) (x-\xo) + A^\top (z - \zo))^\top A^\top u_z
    \\ 
    & + \beta_2 (x - \xo)^\top u_x  + \beta_1 (A(x - \xo))^\top
     A u_x  + \beta_2 (z - \zo)^\top u_z.
  \end{align*}
  Upper bounding now the terms using $\norm{\gradient^2 f(x)},
  \norm{U(\xo,x)} \le M$ for all $x \in \real^n$ yields
  \begin{align}
    & \Lie_{\SDp} V_3(x,z) \notag
    \\
    & \quad \le - [x-\xo; \, A^\top(z-\zo)]^\top \overline{U}(\xo,x)
    [x - \xo; \, A^\top(z - \zo)] \notag
    \\
    & \quad \quad + C_x(x,z) \norm{u_x} + C_z (x,z) \norm{u_z},
    \label{eq:lie-up-bound}
  \end{align}
  where 
  \begin{align*}
    C_x (x,z)& = \Bigl( \beta_1 M^2 \norm{x-\xo} +\beta_1 M \norm{A}
     \norm{z -\zo} 
     \\ 
     & \qquad + \beta_2 \norm{x - \xo} 
         + \beta_1 \norm{A}^2 \norm{x -\xo} \Bigr),
     \\
    C_z (x,z) & = \Bigl( \beta_1 M \norm{A} \norm{x - \xo} +
     \beta_1 \norm{A}^2 \norm{z - \zo}
     \\ 
     & \qquad + \beta_2 \norm{z - \zo}\Bigr),
   \end{align*}
   and $\overline{U}(\xo,x)$ is 
  \begin{align*}
    \begin{bmatrix} \beta_1 U \gradient^2
      f(x) U + \beta_2 U & \beta_1 U \gradient^2 f(x)
      \\
      \beta_1 \gradient^2 f(x) U & \beta_1 \gradient^2
      f(x) \end{bmatrix}.
  \end{align*}
  where $U = U(\xo,x)$.  Note that $C_x(x,z) \le \tilde{C}_x \norm{x -
  \xo; z-\zo} = \tilde{C}_x \norm{(x,z)}_{\sdl{F}}$ and $C_z (x,z) \le
  \tilde{C}_z \norm{x - \xo; z-\zo} = \tilde{C}_z
  \norm{(x,z)}_{\sdl{F}}$, where
  \begin{align*}
    \tilde{C}_x & = \beta_1 M^2 + \beta_1 M \norm{A} + \beta_2 + \beta_1
    \norm{A}^2,
    \\
    \tilde{C}_z & = \beta_1 M \norm{A} + \beta_1 \norm{A}^2 + \beta_2.
  \end{align*}
  From Lemma~\ref{le:U-bound-new}, we have $\overline{U}(\xo,x)
  \succeq \lm_m I$, where $\lm_m > 0$. Employing these facts
  in~\eqref{eq:lie-up-bound}, we obtain
  \begin{align*}
    \Lie_{\SDp} V_3(x,z) 
    & \le - \lm_m (\norm{x-\xo}^2 + \norm{A^\top(z-\zo)}^2)
    \\
    & \qquad + (\tilde{C}_x + \tilde{C}_z)
    \norm{(x,z)}_{\sdl{F}} \norm{u}
  \end{align*}
  From Lemma~\ref{le:A-z-bound}, we get
  \begin{align*}
    \Lie_{\SDp} V_3(x,z)
    & \le - \lm_m(\norm{x-\xo}^2 + \lm_{\mathrm{s}}(A A^\top)
    \norm{z-\zo}^2
    \\
    & \qquad + (\tilde{C}_x + \tilde{C}_z)
    \norm{(x,z)}_{\sdl{F}} \norm{u}
    \\
    & \le - \tilde{\lm}_m \norm{(x,z)}^2_{\sdl{F}}
    \\
    & \qquad + (\tilde{C}_x + \tilde{C}_z)
    \norm{(x,z)}_{\sdl{F}} \norm{u},
  \end{align*}
  where $\tilde{\lm}_m = \lm_m \min\{1, \lm_{\mathrm{s}}(A A^\top)
  \}$. Now pick any $\theta \in (0,1)$. Then,
  \begin{align*}
    \Lie_{\SDp} V_3(x,z) &\le - (1-\theta) \tilde{\lm}_m
    \norm{(x,z)}^2_{\sdl{F}} 
    \\ 
    & \qquad - \theta \tilde{\lm}_m \norm{(x,z)}^2_{\sdl{F}} 
    \\
    & \qquad + (\tilde{C}_x + \tilde{C}_z) \norm{(x,z)}_{\sdl{F}} \norm{u}
    \\
    & \le - (1 - \theta) \tilde{\lm}_m \norm{(x,z)}^2_{\sdl{F}},
  \end{align*}
  whenever $\norm{(x,z)}_{\sdl{F}} \ge \frac{\tilde{C}_x +
  \tilde{C}_z}{\theta \tilde{\lm}_m} \norm{u}$, which proves the ISS
  property.  
\end{IEEEproof}

\begin{remark}\longthmtitle{Relaxing global bounds on Hessian of
    $f$}\label{re:Hess-bound}
  {\rm The assumption on the Hessian of $f$ in
    Theorem~\ref{th:iss-lyapunov} is restrictive, but there are
    functions other than quadratic that satisfy it, see
    e.g.~\cite[Section 6]{SSK-JC-SM:15-auto}. We
    conjecture that the global upper bound on the Hessian can be
    relaxed by resorting to the notion of semiglobal ISS, and we will
    explore this in the future. \oprocend }
\end{remark}

The above result has the following consequence. 

\begin{corollary}\longthmtitle{Lyapunov function for saddle-point
    dynamics}\label{cr:iss-lyapunov-new}
  Let the saddle function $F$ be of the form~\eqref{eq:F-wo-y}, with
  $f$ strongly convex, twice continuously differentiable, and
  satisfying $mI \preceq \gradient^2 f(x) \preceq MI$ for all $x \in
  \real^n$ and some constants $0<m \le M < \infty$.  Then, the
  function $V_3$~\eqref{eq:V2} is a Lyapunov function with respect to
  the set $\sdl{F}$ for the saddle-point
  dynamics~\eqref{eq:saddle-dyn}.
\end{corollary}

\begin{remark}\longthmtitle{ISS with respect to $\sdl{F}$ does not
    imply bounded trajectories}\label{re:unbounded}
  {\rm Note that Theorem~\ref{th:iss-lyapunov} bounds only the
    distance of the trajectories of~\eqref{eq:noise-saddle-dyn}
    to~$\sdl{F}$. Thus, if $\sdl{F}$ is unbounded, the trajectories
    of~\eqref{eq:noise-saddle-dyn} can be unbounded under arbitrarily
    small constant disturbances. However, if matrix $A$ has full
    row-rank, then $\sdl{F}$ is a singleton and the ISS property
    implies that the trajectory of~\eqref{eq:noise-saddle-dyn} remains
    bounded under bounded disturbances.}  \oprocend
\end{remark}

As pointed out in the above remark, if $\sdl{F}$ is not unique, then the
trajectories of the dynamics might not be bounded. We next look at a
particular type of disturbance input which guarantees bounded
trajectories even when $\sdl{F}$ is unbounded. Pick any $(\xo,\zo) \in
\sdl{F}$ and define the function $\map{\tilde{V}_3}{\real^n \times
\real^m}{\realnonnegative}$ as
\begin{multline*}
  \tilde{V}_3(x,z) = \frac{\beta_1}{2}\norm{\SD(x,z)}^2
  + \frac{\beta_2}{2}(\norm{x-\xo}^2 + \norm{z - \zo}^2)
\end{multline*}
with $\beta_1 > 0$, $\beta_2 = \frac{4 \beta_1 M^4}{m^2}$.  One can
show, following similar steps as those of proof of
Theorem~\ref{th:iss-lyapunov}, that the function $\tilde{V}_3$ is an
ISS Lyapunov function with respect to the point $(\xo,\zo)$ for the
dynamics $\SDp$ when the disturbance input to $z$-dynamics has the
special structure $u_z = A \tilde{u}_z$, $\tilde{u}_z \in
\real^n$. This type of disturbance is motivated by scenarios with
measurement errors in the values of $x$ and $z$ used
in~\eqref{eq:saddle-dyn} and without any computation error of the
gradient term in the $z$-dynamics.  The following statement makes
precise the ISS property for this particular disturbance. 

\begin{corollary}\longthmtitle{ISS of saddle-point dynamics}
  \label{cr:iss-lyapunov-2}
  Let the saddle function $F$ be of the form~\eqref{eq:F-wo-y}, with
  $f$ strongly convex, twice continuously differentiable, and
  satisfying $mI \preceq \gradient^2 f(x) \preceq MI$ for all $x \in
  \real^n$ and some constants $0<m \le M < \infty$.  Then, the
  dynamics
  \begin{align}\label{eq:noise-saddle-dyn-2}
    \begin{bmatrix}
      \dot x
      \\
      \dot z
    \end{bmatrix}
    =
    \begin{bmatrix}
      - \gradient_x F(x,z)
      \\
      \gradient_z F(x,z)
    \end{bmatrix}
    +
    \begin{bmatrix}
      u_x
      \\
      A \tilde{u}_z
    \end{bmatrix},
  \end{align}
  where $\map{(u_x, \tilde{u}_z)}{\realnonnegative}{\real^{2n}}$ is
  measurable and locally essentially bounded input, is ISS with
  respect to every point of $\sdl{F}$.
\end{corollary}

The proof is analogous to that of Theorem~\ref{th:iss-lyapunov} with the
key difference that the terms $C_x(x,z)$ and $C_z(x,z)$ appearing
in~\eqref{eq:lie-up-bound} need to be upper bounded in terms of
$\norm{x-\xo}$ and $\norm{A^\top(z-\zo)}$. This can be done due to the
special structure of $u_z$. With these bounds, one arrives at the
condition~\eqref{eq:iss-lyap-c2} for Lyapunov $\tilde{V}_3$ and
dynamics~\eqref{eq:noise-saddle-dyn-2}. One can deduce from
Corollary~\ref{cr:iss-lyapunov-2} that the trajectory
of~\eqref{eq:noise-saddle-dyn-2} remains bounded for bounded input even
when $\sdl{F}$ is unbounded.

\changes{
\begin{example}\longthmtitle{ISS property of saddle-point dynamics}
  {\rm
  Consider $\map{F}{\real^2 \times \real^2}{\real}$ of the
  form~\eqref{eq:F-wo-y} with
  \begin{align}\label{eq:f-A-b}
    f(x) = x_1^2 + (x_2 - 2)^2, \notag
    \\
    A =
    \begin{bmatrix}
      1 & -1
      \\
      -1 & 1 
    \end{bmatrix}, \text{ and } b =
    \begin{bmatrix}
      0
      \\
      0
    \end{bmatrix}.
  \end{align}
  Then, $\sdl{F} = \setdef{(x,z) \in \real^2 \times \real^2}{x =
    (1,1), z =  (0,2) + \lambda (1,1), \lambda \in \real }$ is
  a continuum of points. Note that $\gradient^2 f(x) = 2I$, thus,
  satisfying the assumption of bounds on the Hessian of $f$. By
  Theorem~\ref{th:iss-lyapunov}, the saddle-point dynamics for this
  saddle function $F$ is input-to-state stable with respect to the set
  $\sdl{F}$.  This fact is illustrated in Figure~\ref{fig:iss-two},
  which also depicts how the specific structure of the disturbance
  input in~\eqref{eq:noise-saddle-dyn-2} affects the boundedness of
  the trajectories. 
  \begin{figure}[htb]
    \centering%
    \subfloat[$(x,z)$]{\includegraphics[width = 0.475
      \linewidth]{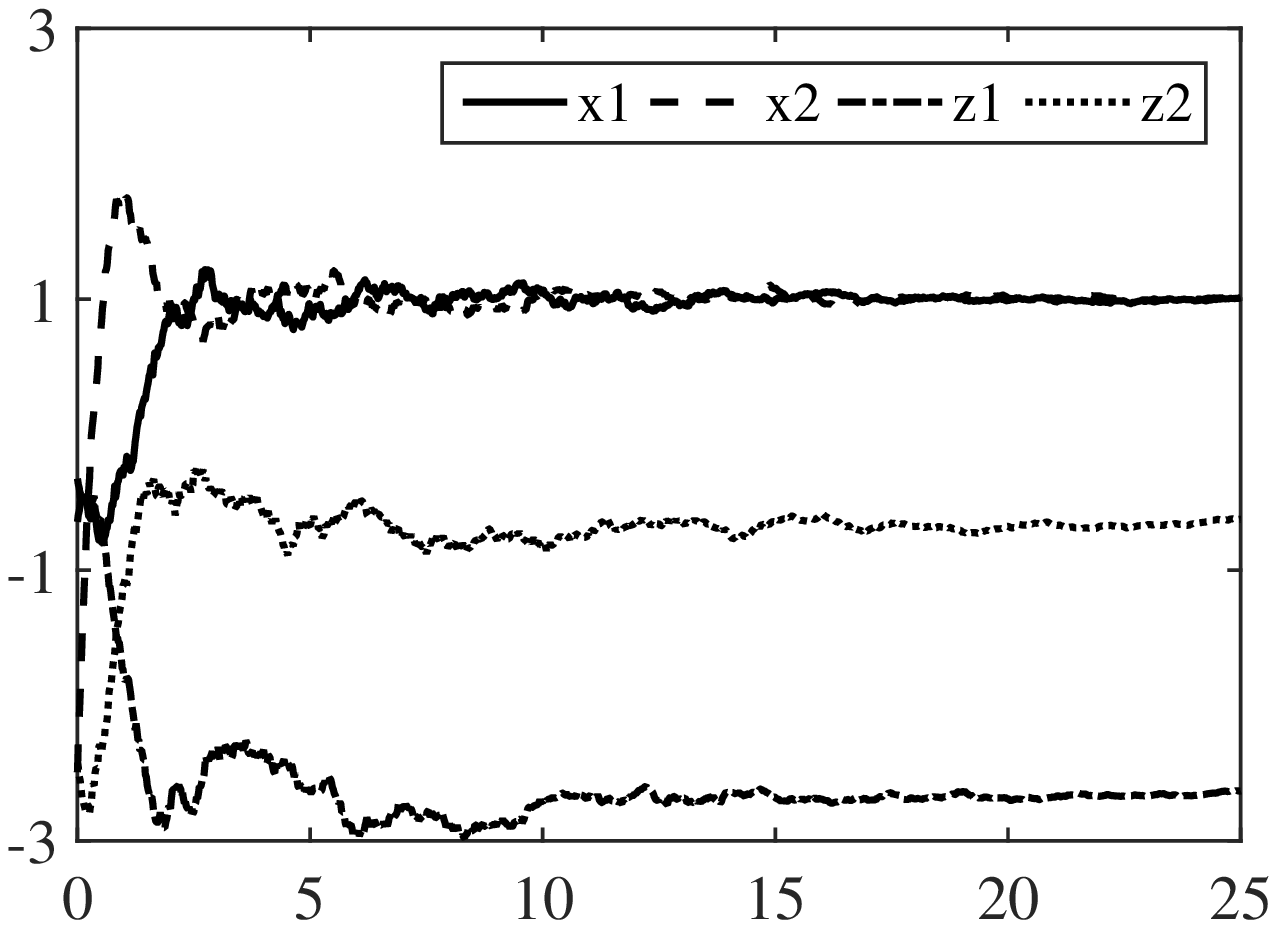}}
    \quad
    \subfloat[$\norm{(x,z)}_{\sdl{F}}$]{\includegraphics[width = 0.475
      \linewidth]{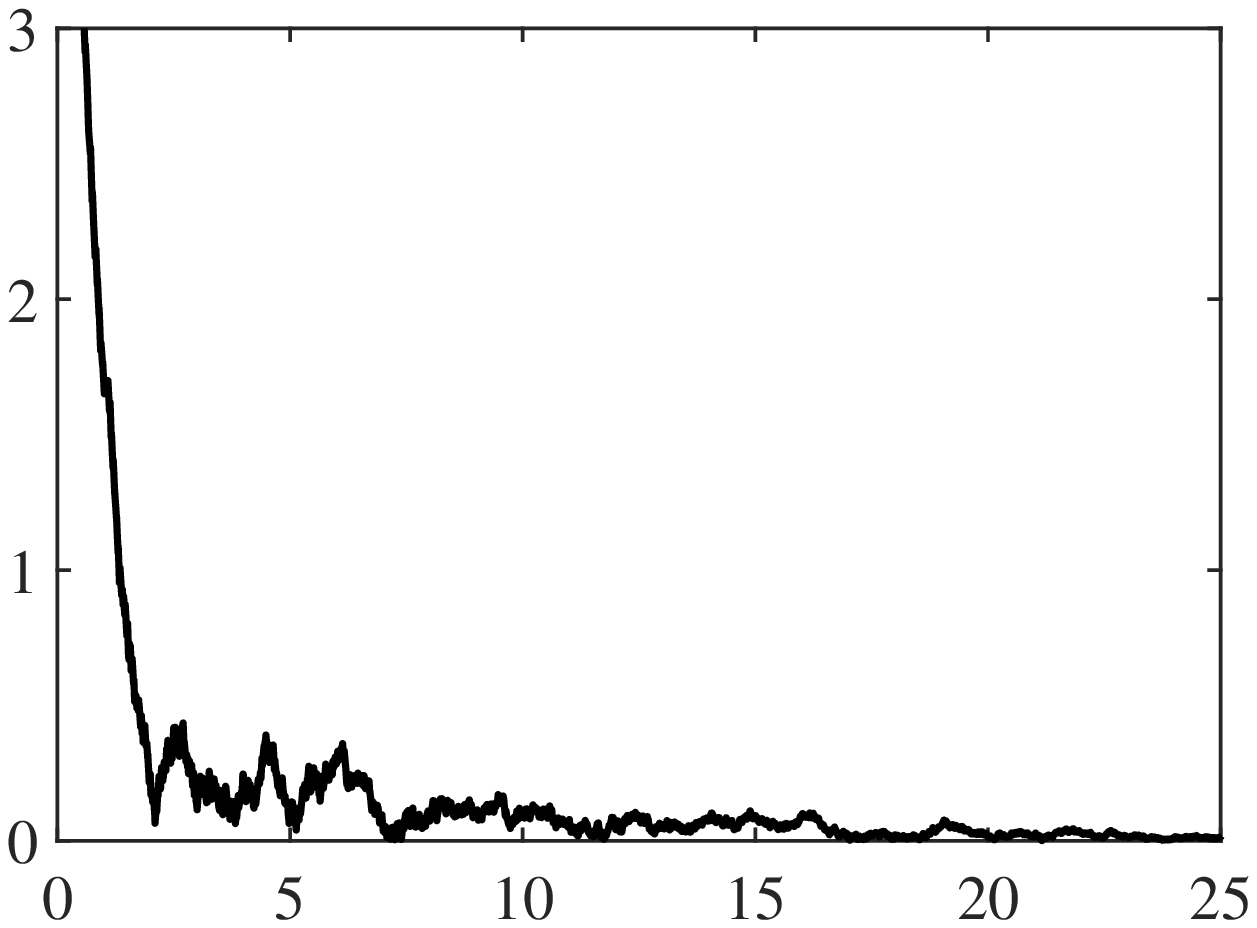}}
    \\
    \subfloat[$(x,z)$]{\includegraphics[width = 0.475
      \linewidth]{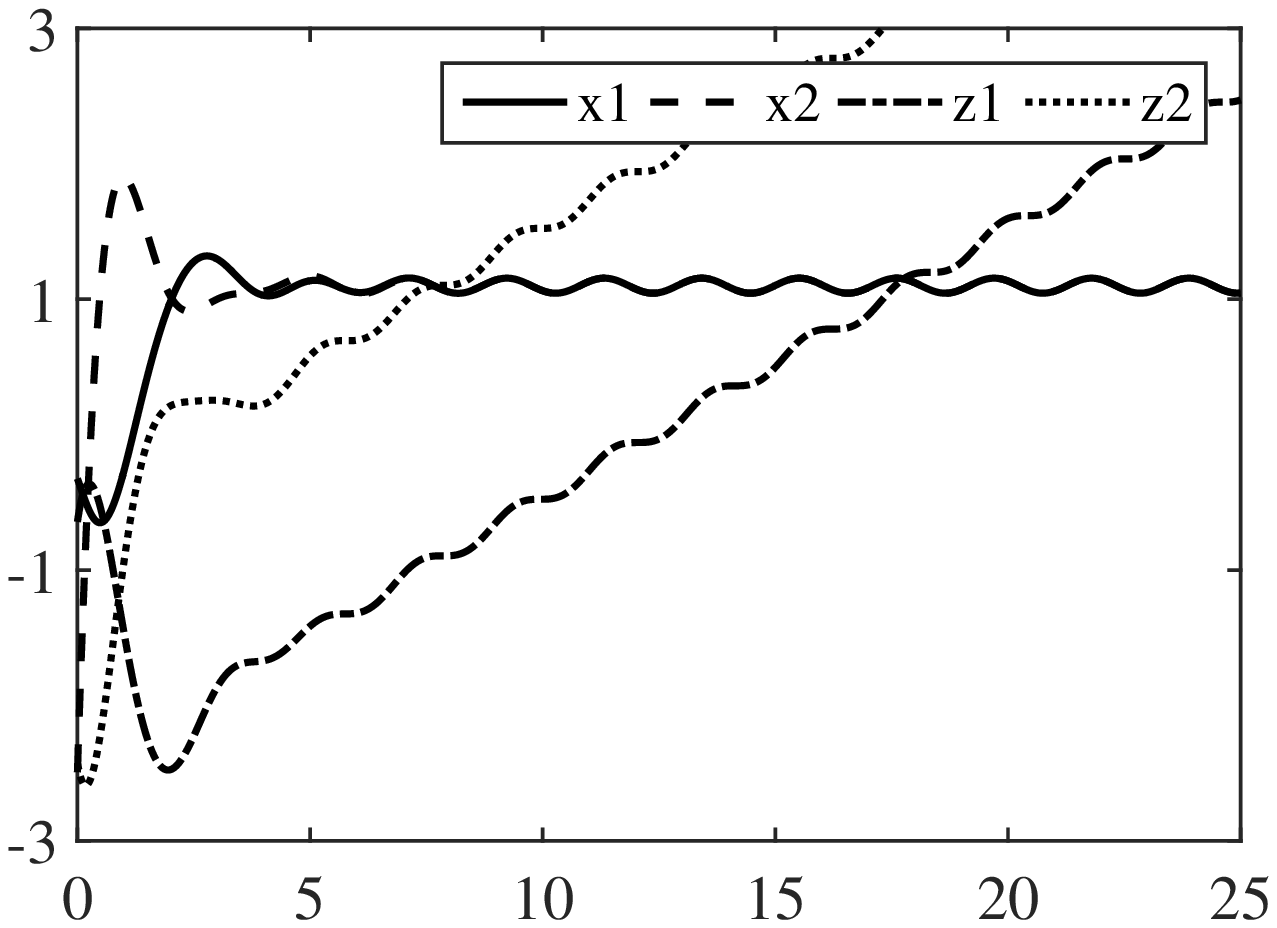}}
    \quad
    \subfloat[$\norm{(x,z)}_{\sdl{F}}$]{\includegraphics[width = 0.475
      \linewidth]{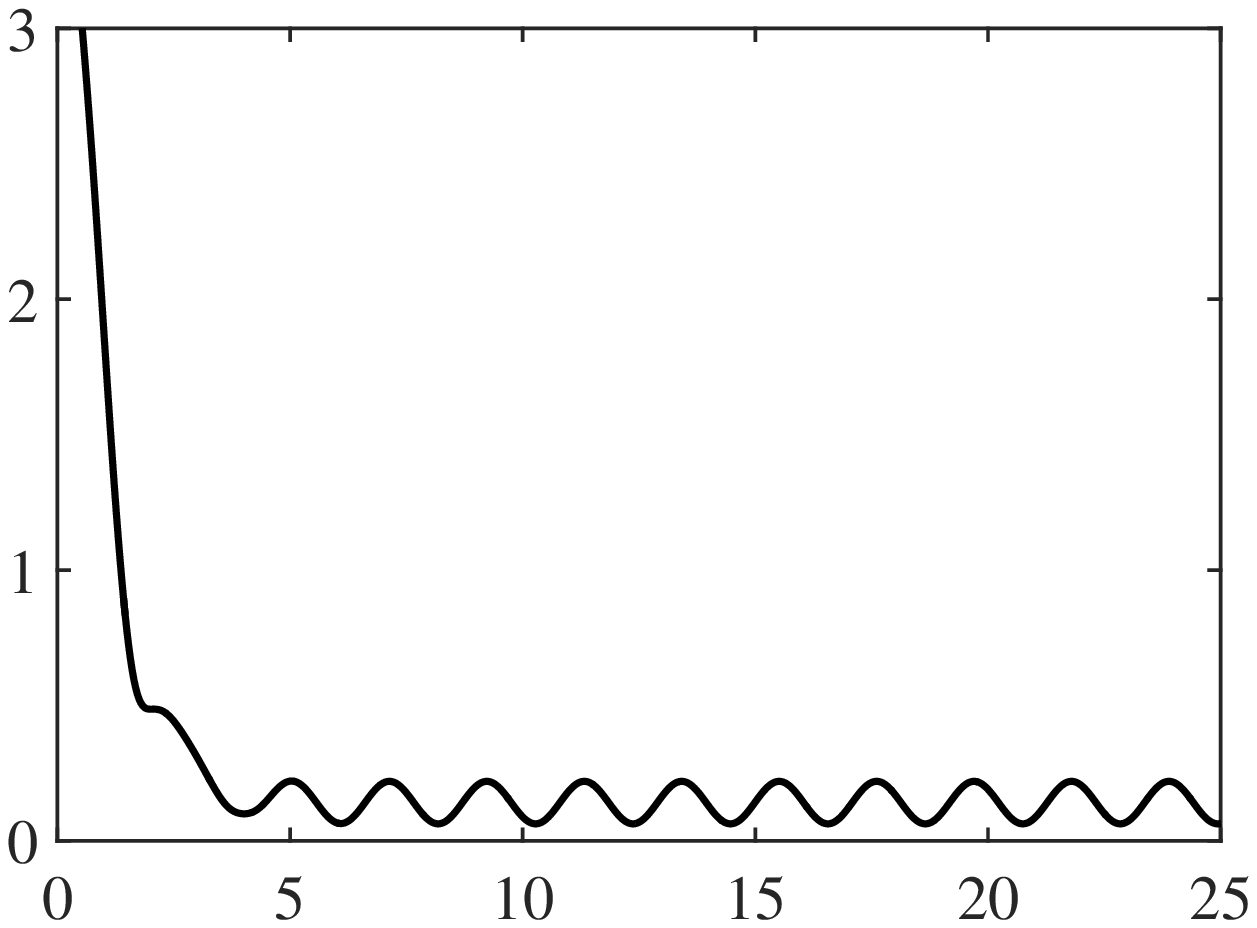}}
    \\
    \subfloat[$(x,z)$]{\includegraphics[width = 0.475
      \linewidth]{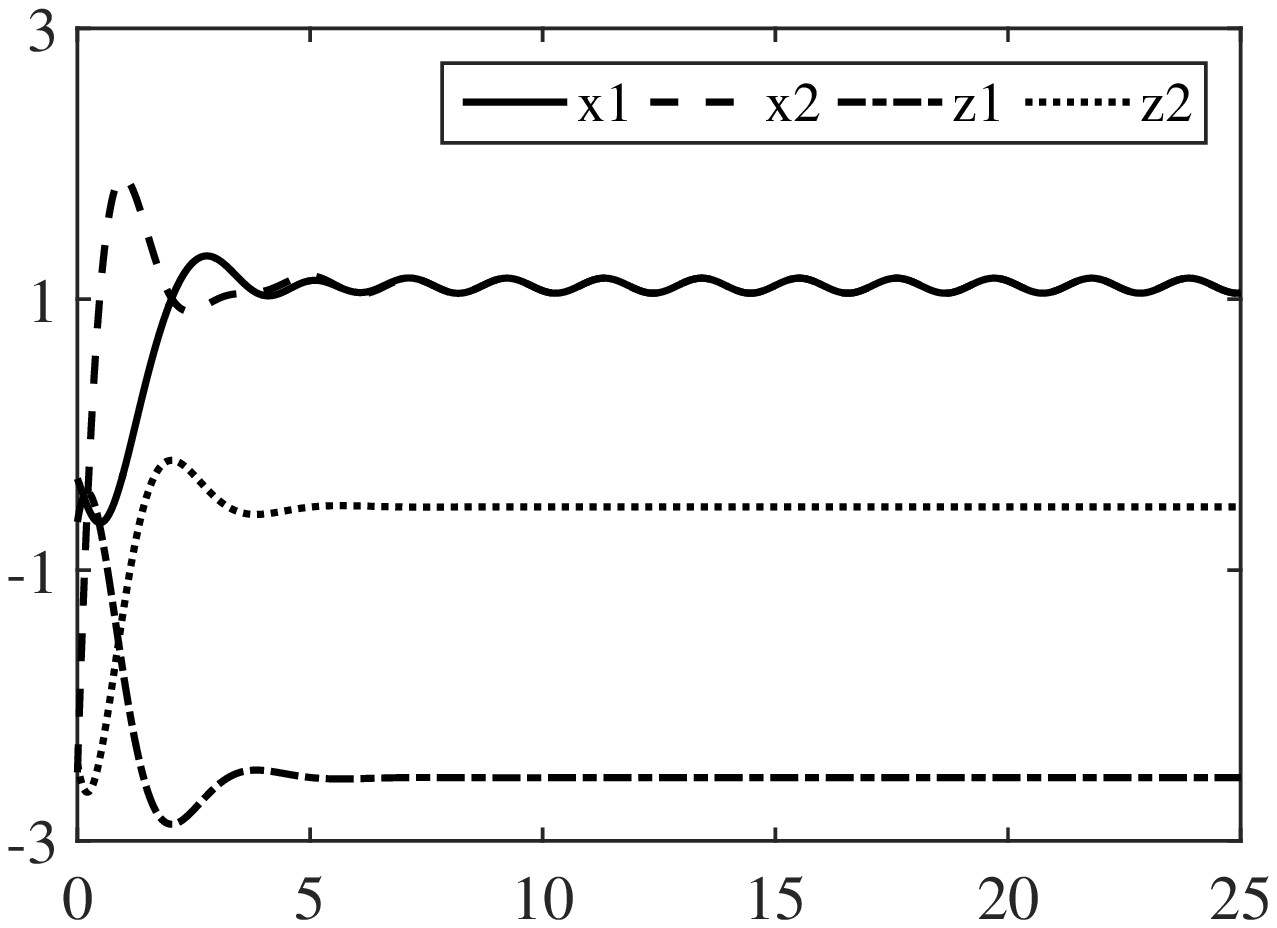}}
    \quad
    \subfloat[$\norm{(x,z)}_{\sdl{F}}$]{\includegraphics[width = 0.475
      \linewidth]{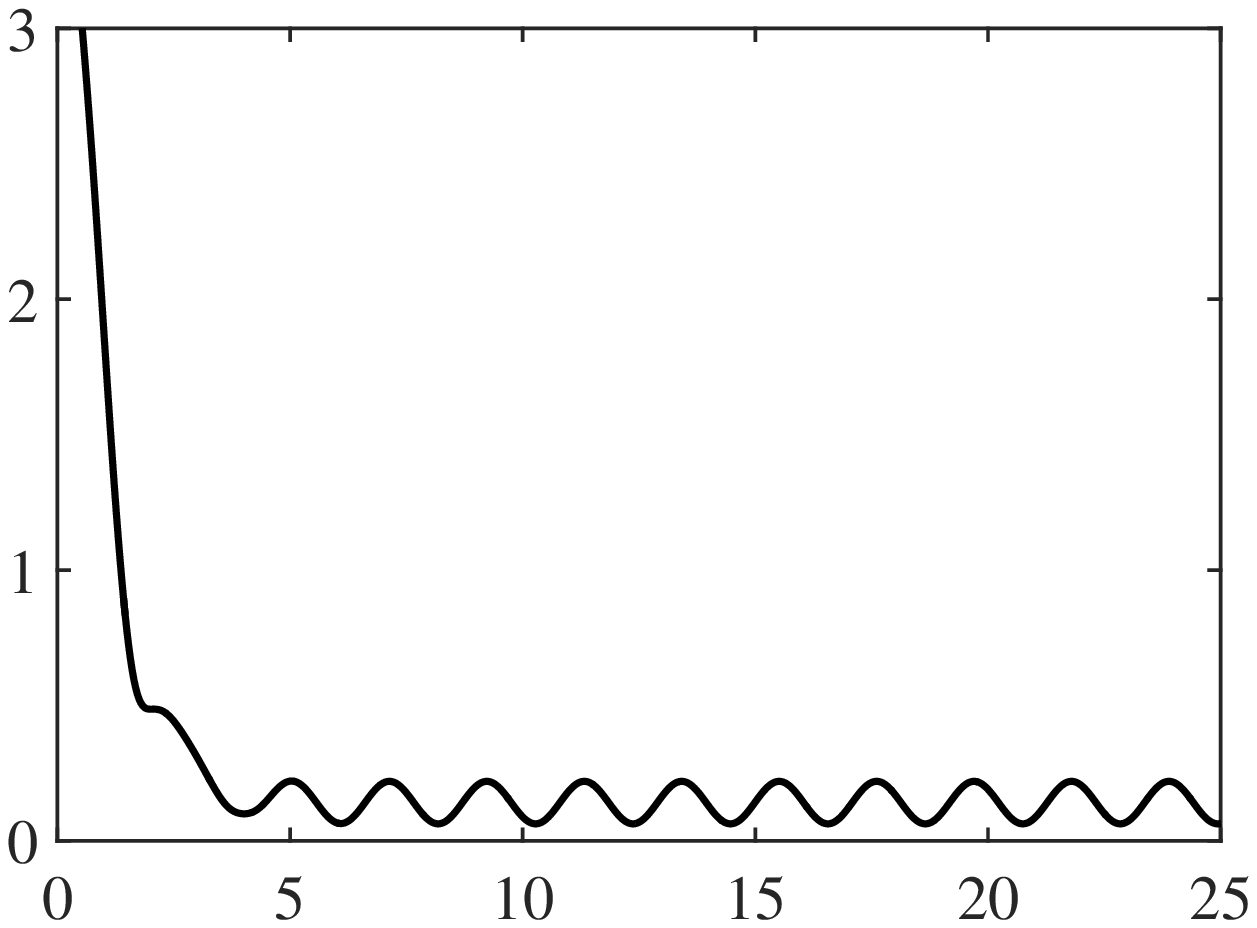}}
    \caption{Plots (a)-(b) show the ISS property, cf
      Theorem~\ref{th:iss-lyapunov}, of the
      dynamics~\eqref{eq:noise-saddle-dyn} for the saddle function $F$
      defined by~\eqref{eq:f-A-b}. The initial condition is $x(0) =
      (-0.3254, -2.4925)$ and $z(0) = (-0.6435, -2.4234)$ and the
      input $u$ is exponentially decaying in magnitude. As shown in
      (a)-(b), the trajectory converges asymptotically to a saddle
      point as the input is vanishing.  Plots (c)-(d) have the same
      initial condition but the disturbance input consists of a
      constant plus a sinusoid.  The trajectory is unbounded under
      bounded input while the distance to the set of saddle points
      remains bounded, cf. Remark~\ref{re:unbounded}.  Plots (e)-(f)
      have the same initial condition but the disturbance input to the
      $z$-dynamics is of the form~\eqref{eq:noise-saddle-dyn-2}. In
      this case, the trajectory remains bounded as the dynamics is ISS
      with respect to each saddle point, cf.
      Corollary~\ref{cr:iss-lyapunov-2}.  }\label{fig:iss-two}
  \end{figure}
  \oprocend }
\end{example}
}
\changes{
  \begin{remark}\longthmtitle{Quadratic ISS-Lyapunov
      function}\label{re:alternative-ISS}
  {\rm For the saddle-point dynamics~\eqref{eq:saddle-dyn}, the ISS
    property stated in Theorem~\ref{th:iss-lyapunov} and
    Corollary~\ref{cr:iss-lyapunov-2} can also be shown using a
    quadratic Lyapunov function. Let
    $\map{V_4}{\real^n \times \real^m}{\realnonnegative}$ be
    \begin{align*}
      V_4(x,z) = \frac{1}{2} \norm{(x,z)}_{\sdl{F}}^2 + \eps (x -
      x_p)^\top A^\top (z-z_p), 
    \end{align*}
    where $(x_p,z_p) = \proj_{\sdl{F}}(x,z)$ and $\eps > 0$. Then, one
    can show that there exists $\eps_{\max} > 0$ such that $V_4$ for
    any $\eps \in (0,\eps_{\max})$ is an ISS-Lyapunov function for the
    dynamics~\eqref{eq:saddle-dyn}. For space reasons, we omit the
    complete analysis of this fact here.  \oprocend }
\end{remark}
}
\subsection{Self-triggered implementation}

In this section we develop an opportunistic state-triggered
implementation of the (continuous-time) saddle-point dynamics. Our aim
is to provide a discrete-time execution of the algorithm, either on a
physical system or as an optimization strategy, that do not require
the continuous evaluation of the vector field and instead adjust the
stepsize based on the current state of the system. Formally, given a
sequence of triggering time instants $\{t_k\}_{k=0}^\infty$, with $t_0
= 0$, we consider the following implementation of the saddle-point
dynamics
\begin{subequations}\label{eq:self-trig-dyn}
  \begin{align}
    \dot x(t) & = - \gradient_x F(x(t_k),z(t_k)),
    \\
    \dot z(t) & = \gradient_z F(x(t_k),z(t_k)).
  \end{align}
\end{subequations}
for $t \in [t_k, t_{k+1})$ and $k \in \integersnonnegative$. The
objective is then to design a criterium to opportunistically select
the sequence of triggering instants, guaranteeing at the same time the
feasibility of the execution and global asymptotic convergence, see
e.g.,~\cite{WPMHH-KHJ-PT:12}. Towards this goal, we look at the
evolution of the Lyapunov function~$V_3$ in~\eqref{eq:V2}
along~\eqref{eq:self-trig-dyn},
\begin{align}\label{eq:self-trig-design}
  & \gradient V_3(x(t),z(t))^\top \SD(x(t_k),z(t_k)) \notag
  \\
  & = \Lie_{\SD} V_3(x(t_k),z(t_k))
  \\
  & \quad + \Bigl(\gradient V_3(x(t),z(t)) - \gradient
  V_3(x(t_k),z(t_k)) \Bigr)^\top \notag
  \\
  & \qquad \qquad \qquad \qquad \qquad \qquad \qquad \qquad \quad
  \SD(x(t_k),z(t_k)) . \notag
\end{align}
We know from Corollary~\ref{cr:iss-lyapunov-new} that the first summand
is negative outside $\sdl{F}$. Clearly, for $t=t_k$, the second
summand vanishes, and by continuity, for $t$ sufficiently close to
$t_{k}$, this summand remains smaller in magnitude than the first,
ensuring the decrease of~$V_3$. \changes{To make this argument precise,
  we employ Proposition~\ref{prop:V-bound}
  in~\eqref{eq:self-trig-design} and obtain}
\begin{align*}
  & \gradient V_3(x(t),z(t))^\top
  \SD(x(t_k),z(t_k))
  \\
  & \le  \Lie_{\SD} V_3(x(t_k),z(t_k))+ \xi(x(t_k),z(t_k)) 
  \\
  & \quad    \norm{(x(t) - x(t_k)); (z(t) - z(t_k))}
  \norm{\SD(x(t_k),z(t_k))}
  \\
  & = \Lie_{\SD} V_3(x(t_k),z(t_k)) 
  \\
  & \qquad + (t - t_k) \xi(x(t_k),z(t_k))
  \norm{\SD(x(t_k),z(t_k))}^2,
\end{align*}
where the equality follows from writing $(x(t),z(t))$ in terms of
$(x(t_k),z(t_k))$ by integrating~\eqref{eq:self-trig-dyn}.  Therefore,
in order to ensure the monotonic decrease of $V_3$, we require the above
expression to be nonpositive. That is, 
\begin{align}\label{eq:trig-time-pre}
  t_{k+1} \le t_k - \frac{\Lie_{\SD} V_3(x(\tht_k),z(\tht_k))}
  {\xi(x(\tht_k),z(\tht_k)) \norm{\SD(x(\tht_k),z(\tht_k))}^2}.
\end{align}
Note that to set $t_{k+1}$ equal to the right-hand side of the above
expression, one needs to compute the Lie derivative
at~$(x(t_k),z(t_k))$. We then distinguish between two
possibilities. If the self-triggered saddle-point dynamics acts as a
closed-loop physical system and its equilibrium points are known, then
computing the Lie derivative is feasible and one can
use~\eqref{eq:trig-time-pre} to determine the triggering times. If,
however, the dynamics is employed to seek the primal-dual optimizers
of an optimization problem, then computing the Lie derivative is
infeasible as it requires knowledge of the optimizer. To overcome this
limitation, we propose the following alternative triggering criterium
which satisfies~\eqref{eq:trig-time-pre} as shown later in our
convergence analysis,
\begin{align}\label{eq:trig-time}
  t_{k+1} = t_k +  \frac{\tilde{\lm}_m}{3(M^2 + \norm{A}^2)
\xi(x(t_k),z(t_k))},
\end{align}
where $\tilde{\lm}_m = \lm_m \min\{1,\lm_{\mathrm{s}}(A A^\top)\}$,
$\lm_m$ is given in Lemma~\ref{le:U-bound-new}, and $\lm_{\mathrm{s}}
(A A^\top)$ is the smallest nonzero eigenvalue of $A A^\top$.  In
either~\eqref{eq:trig-time-pre} or~\eqref{eq:trig-time}, the
right-hand side depends only on the state $(x(t_k),z(t_k))$. These
triggering times for the dynamics~\eqref{eq:self-trig-dyn} define a
first-order Euler discretization of the saddle-point dynamics with
step-size selection based on the current state of the system.
It is for this reason that we refer to~\eqref{eq:self-trig-dyn}
together with either the triggering criterium~\eqref{eq:trig-time-pre}
or~\eqref{eq:trig-time} as the \emph{self-triggered saddle-point
  dynamics}.  In integral form, this dynamics results in a
discrete-time implementation of~\eqref{eq:saddle-dyn} given as
\begin{align*}
  \begin{bmatrix} x(t_{k+1})
    \\
    z(t_{k+1})
  \end{bmatrix}
  = 
  \begin{bmatrix}
    x(t_k)
    \\
    z(t_k)
  \end{bmatrix}
  + (t_{k+1} - t_k) \SD(x(t_k),z(t_k)) .
\end{align*}
\changes{Note that this dynamics can also be regarded as a
  state-dependent switched system with a single continuous mode and a
  reset map that updates the sampled state at the switching times,
  cf.~\cite{DL:03}.}
We understand the solution of~\eqref{eq:self-trig-dyn} in the
Caratheodory sense (note that this dynamics has a discontinuous
right-hand side). The existence of such solutions, possibly defined only
on a finite time interval, is guaranteed from the fact that along any
trajectory of the dynamics there are only countable number of
discontinuities encountered in the vector field.  The next result
however shows that solutions of~\eqref{eq:self-trig-dyn} exist over the
entire domain $[0,\infty)$ as the difference between consecutive
  triggering times of the solution is lower bounded by a positive
  constant. Also, it establishes the asymptotic convergence of solutions
  to the set of saddle points. 

\begin{theorem}\longthmtitle{Convergence of the self-triggered
    saddle-point dynamics}\label{th:conv-self-trig-dyn}
  Let the saddle function $F$ be of the form~\eqref{eq:F-wo-y}, with $A$
  having full row rank, $f$ strongly convex, twice 
  differentiable, and satisfying $m I \preceq \gradient^2 f(x) \preceq M
  I$ for all $x \in \real^n$ and some constants $0< m \le M < \infty$.
  Let the map $x \mapsto \gradient^2 f(x)$ be Lipschitz
  with some constant $L>0$. Then, $\sdl{F}$ is singleton. Let $\sdl{F}
  = \{(\xo,\zo)\}$. Then, for any initial condition $(x(0),z(0))
  \in \real^n \times \real^m$, we have
  \begin{align*}
    \lim_{k \to \infty} (x(t_k),z(t_k)) = (\xo,\zo)
  \end{align*}
  for the solution of the self-triggered saddle-point dynamics,
  defined by~\eqref{eq:self-trig-dyn} and~\eqref{eq:trig-time}, starting at
  $(x(0),z(0))$. Further, there exists $\mu_{(x(0),z(0))} > 0$ such that
  the triggering times of this solution satisfy 
   \begin{align*}
    t_{k+1} - t_k \ge \mu_{(x(0),z(0))}, \quad \text{ for all } k \in
    \naturalnumbers.
  \end{align*}
\end{theorem}
\begin{IEEEproof}
  Note that there is a unique equilibrium point to the saddle-point
  dynamics~\eqref{eq:saddle-dyn} for $F$ satisfying the stated
  hypotheses. Therefore, the set of saddle point is singleton for this
  $F$. Now, given $(x(0),z(0)) \in \real^n \times \real^m$, let $V_3^0
  = V_3(x(0),z(0))$ and define
  \begin{align*}
    G = \max\setdef{\norm{\gradient_x F(x,z)}}{(x,z) \in
      \levelset{V_3}{V_3^0}},
  \end{align*}
  where, we use the notation for the sublevel set of $V_3$ as 
  \begin{align*}
    \levelset{V_3}{\alpha} =\setdef{(x,z) \in \real^n \times
    \real^m}{V_3(x,z) \le \alpha}
  \end{align*}
  for any $\alpha \ge 0$. Since $V_3$ is radially unbounded, the set
  $\levelset{V_3}{V_3^0}$ is compact and so, $G$ is well-defined and
  finite. If the trajectory of the self-triggered saddle-point dynamics
  is contained in $\levelset{V_3}{V_3^0}$, then we can bound the
  difference between triggering times in the following way. 
  From Proposition~\ref{prop:V-bound} for all $(x,z) \in
  \levelset{V_3}{V_3^0}$, we have $\xi_1(x,z) = M \xi_2 + L \norm{
  \gradient_x F(x,z)} \le M \xi_2 + L G =: T_1$. Hence, for all $(x,z) \in
  \levelset{V_3}{V_3^0}$, we get 
  \begin{align*}
    \xi(x,z) & = \Bigl( \beta_1^2 (\xi_1(x,z)^2 + \norm{A}^4 + \norm{A}^2
    \xi_2^2) + \beta_2^2 \Bigr)^{\frac{1}{2}}
    \\
    & \le \Bigl( \beta_1^2 (T_1^2 + \norm{A}^4 + \norm{A}^2 + \xi_2^2) +
    \beta_2^2 \Bigr)^{\frac{1}{2}}
    \\
    & =: T_2.
  \end{align*}
  Using the above bound in~\eqref{eq:trig-time}, we get for all $k \in
  \naturalnumbers$
  \begin{align*}
    t_{k+1} - t_k & = \frac{\tilde{\lm}_m}{3(M^2 + \norm{A}^2)
    \xi(x(t_k),z(t_k))}
    \\
    & \ge \frac{\tilde{\lm}_m}{3 (M^2 +\norm{A}^2) T_2} > 0.
  \end{align*}
  This implies that as long as the trajectory is contained in
  $\levelset{V_3}{V_3^0}$, the inter-trigger times are lower bounded by
  a positive quantity. Our next step is to show that the trajectory is
  contained in $\levelset{V_3}{V_3^0}$. Note that
  if~\eqref{eq:trig-time-pre} is satisfied for the triggering
  condition~\eqref{eq:trig-time}, then the sequence
  $\{V_3(x(t_k),z(t_k))\}_{k \in \naturalnumbers}$ is strictly
  decreasing. Since $V_3$ is nonnegative, this implies that $\lim_{k \to
  \infty} V_3(x(t_k),z(t_k)) = 0$ and so, by continuity, $\lim_{k \to \infty}
  (x(t_k),z(t_k)) = (\xo,\zo)$. Thus, it remains to show
  that~\eqref{eq:trig-time} implies~\eqref{eq:trig-time-pre}. To this
  end, first note the following inequalities shown in the proof of
  Theorem~\ref{th:iss-lyapunov}
  \begin{subequations}
    \begin{align}
      \frac{\norm{\SD(x,z)}^2}{3 (M^2 +
	\norm{A}^2)} & \le \norm{(x-\xo);(z-\zo)}^2, \label{eq:bound1} 
	\\
	\abs{\Lie_{\SD} V_3(x,z)} & \ge \tilde{\lm}_m
      \norm{(x-\xo);(z-\zo)}^2. \label{eq:bound2}
    \end{align}
  \end{subequations}
  Using these bounds, we get from~\eqref{eq:trig-time}
  \begin{align*}
    & t_{k+1} - t_k 
    \\
    & = \frac{\tilde{\lm}_m}{3(M^2 + \norm{A}^2) \xi(x(t_k),z(t_k))}
    \\
    & \overset{(a)}{=} \frac{\tilde{\lm}_m
    \norm{\SD(x(t_k),z(t_k))}^2}{3(M^2 + \norm{A}^2) \xi(x(t_k),z(t_k))
    \norm{\SD(x(t_k),z(t_k))}^2}
    \\
    & \overset{(b)}{\le} \frac{\tilde{\lm}_m
    \norm{(x(t_k)-\xo);(z(t_k)-\zo)}^2}{\xi(x(t_k),z(t_k))
      \norm{\SD(x(t_k),z(t_k))}^2}
    \\
    & \overset{(c)}{\le}
    \frac{\abs{\Lie_{\SD}V_3(x(t_k),z(t_k))}}{\xi(x(t_k),z(t_k))
      \norm{\SD(x(t_k),z(t_k))}^2}
    \\
    & = - \frac{\Lie_{\SD} V_3(x(t_k),z(t_k))}{\xi(x(t_k),z(t_k))
    \norm{\SD(x(t_k),z(t_k))}^2},
  \end{align*}
  where (a) is valid as $\norm{\SD(x(t_k),z(t_k))} \not = 0$,
  (b) follows from~\eqref{eq:bound1}, and (c)
  follows from~\eqref{eq:bound2}. Thus,~\eqref{eq:trig-time}
  implies~\eqref{eq:trig-time-pre} which completes the proof. 
\end{IEEEproof}

Note from the above proof that the convergence implication of
Theorem~\ref{th:conv-self-trig-dyn} is also valid when the triggering
criterium is given by~\eqref{eq:trig-time-pre} with the inequality
replaced by the equality.

\begin{example}\longthmtitle{Self-triggered saddle-point
  dynamics}\label{ex:self-triggered}
  {\rm
  Consider the function $\map{F}{\real^3 \times \real}{\real}$, 
  \begin{equation}\label{eq:F-self-trig}
    F(x,z) = \norm{x}^2 + z(x_1 + x_2 + x_3 -1).
  \end{equation}
  Then, with the notation of~\eqref{eq:F-wo-y}, we have $f(x) =
  \norm{x}^2$, $A = [1, 1, 1]$, and $b = 1$. The set of saddle points is
  a singleton, $\sdl{F} = \{
  ((\frac{1}{3},\frac{1}{3},\frac{1}{3}) , -\frac{2}{3}) \}$. Note
  that $\gradient^2 f(x) = 2I$ and $A$ has full row-rank, thus, the
  hypotheses of Theorem~\ref{th:conv-self-trig-dyn} are met. Hence,
  for this $F$, the self-triggered saddle-point
  dynamics~\eqref{eq:self-trig-dyn} with triggering
  times~\eqref{eq:trig-time} converges asymptotically to the saddle
  point of $F$. Moreover, the difference between two consecutive
  triggering times is lower bounded by a finite quantity.
  Figure~\ref{fig:self-triggered} illustrates a simulation of
  dynamics~\eqref{eq:self-trig-dyn} with triggering
  criteria~\eqref{eq:trig-time-pre} (replacing inequality with
  equality), showing that this triggering criteria also ensures
  convergence as commented above. 
  \begin{figure}
    \centering
    \subfloat[$(x,z)$]{\includegraphics[width = 0.475
    \linewidth]{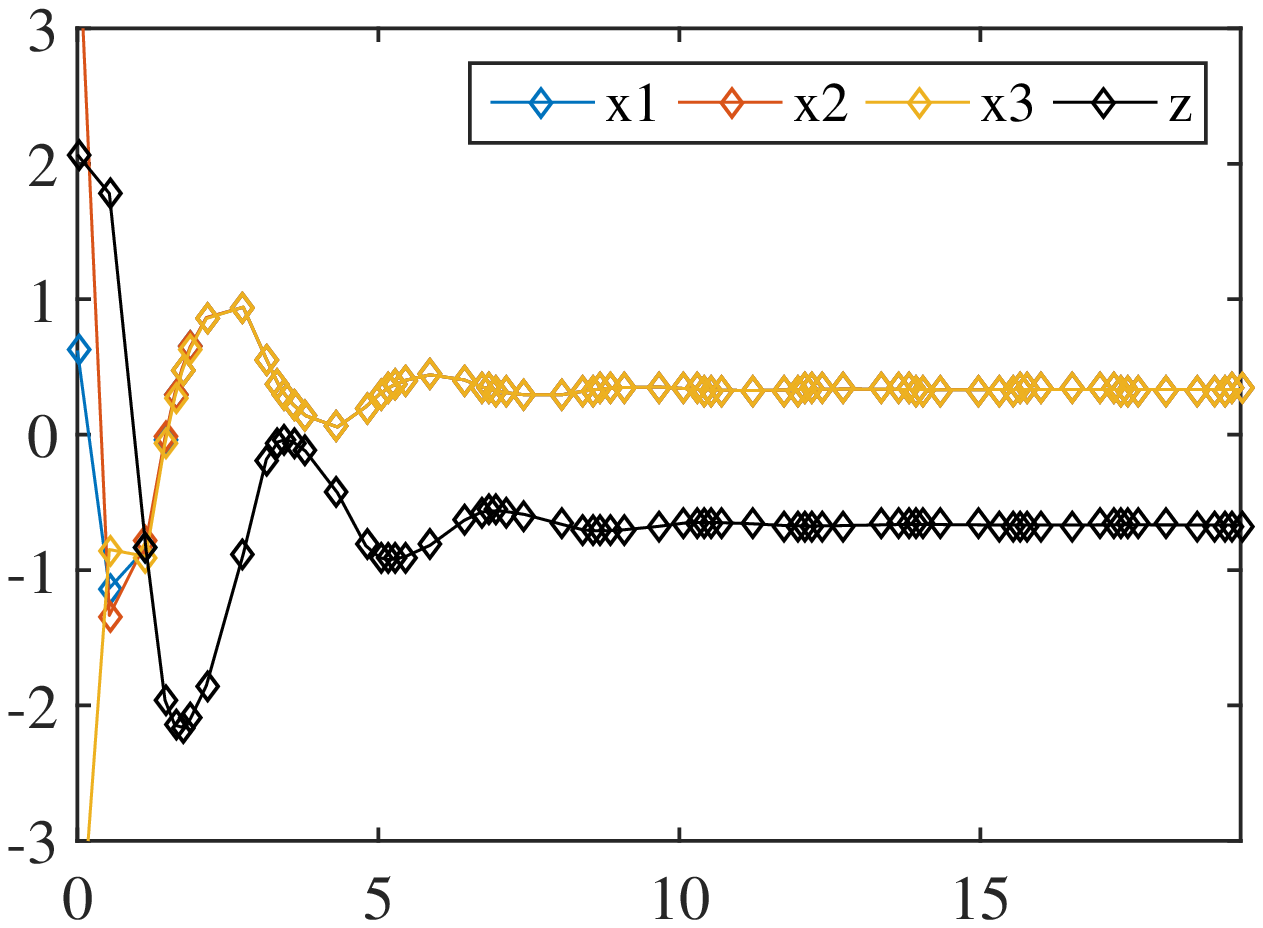}}
    \quad
    \subfloat[$V_3$]{\includegraphics[width = 0.475
    \linewidth]{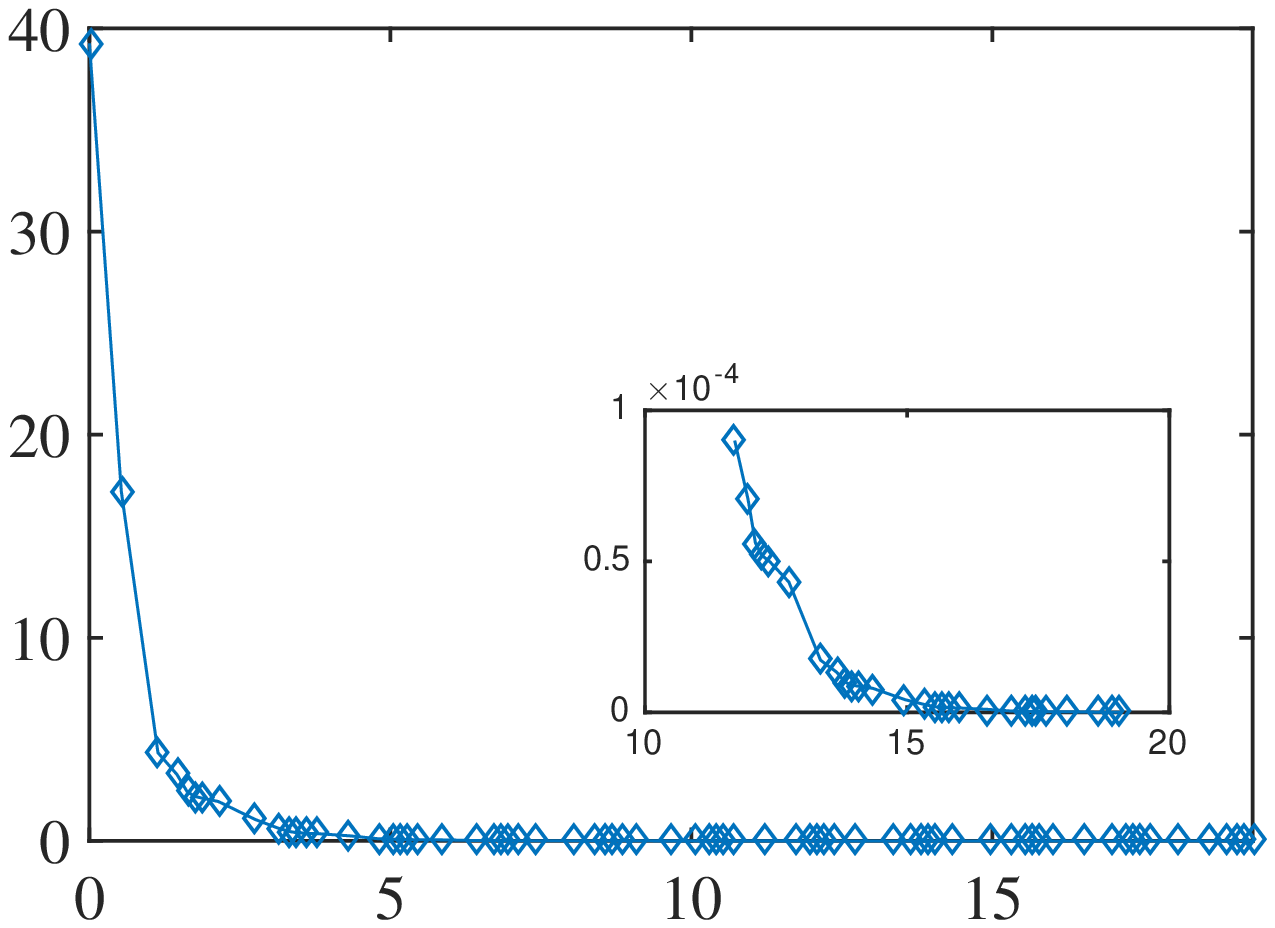}}
  \caption{Illustration of the self-triggered saddle-point dynamics
    defined by~\eqref{eq:self-trig-dyn} with the triggering
    criterium~\eqref{eq:trig-time-pre}. The saddle function $F$ is
    defined in~\eqref{eq:F-self-trig}. With respect to the notation of
    Theorem~\ref{th:conv-self-trig-dyn}, we have $m = M =2$ and
    $\norm{A} = \sqrt{3}$. We select $\beta_1 = 0.1$, then
    $\beta_2 = 1.6$, and from~\eqref{eq:xi-def}, $\xi_1 = 2$. These
    constants define functions $V_3$ (cf.~\eqref{eq:V2}), $\xi$, and
    $\xi_2$ (cf.~\eqref{eq:xi-def}) and also, the triggering
    times~\eqref{eq:trig-time}. In plot (a), the initial condition is
    $x(0) = (0.6210, 3.9201, -4.0817)$, $z(0) = 2.0675$. The
    trajectory converges to the unique saddle-point and the
    inter-trigger times are lower bounded by a positive
    quantity.}\label{fig:self-triggered}
  \end{figure}
  \changes{Finally, Figure~\ref{fig:comparison} compares the self-triggered
  implementation of the saddle-point dynamics with a constant-stepsize
  and a decaying-stepsize first-order Euler discretization. In both
  cases, the the self-triggered dynamics achieves convergence faster,
  and this may be attributed to the fact that it tunes the stepsize in a
  state-dependent way.}
  \begin{figure}
    \centering
    \subfloat[$\norm{(x,z)}_{\sdl{F}}$]{\includegraphics[width = 0.475
      \linewidth]{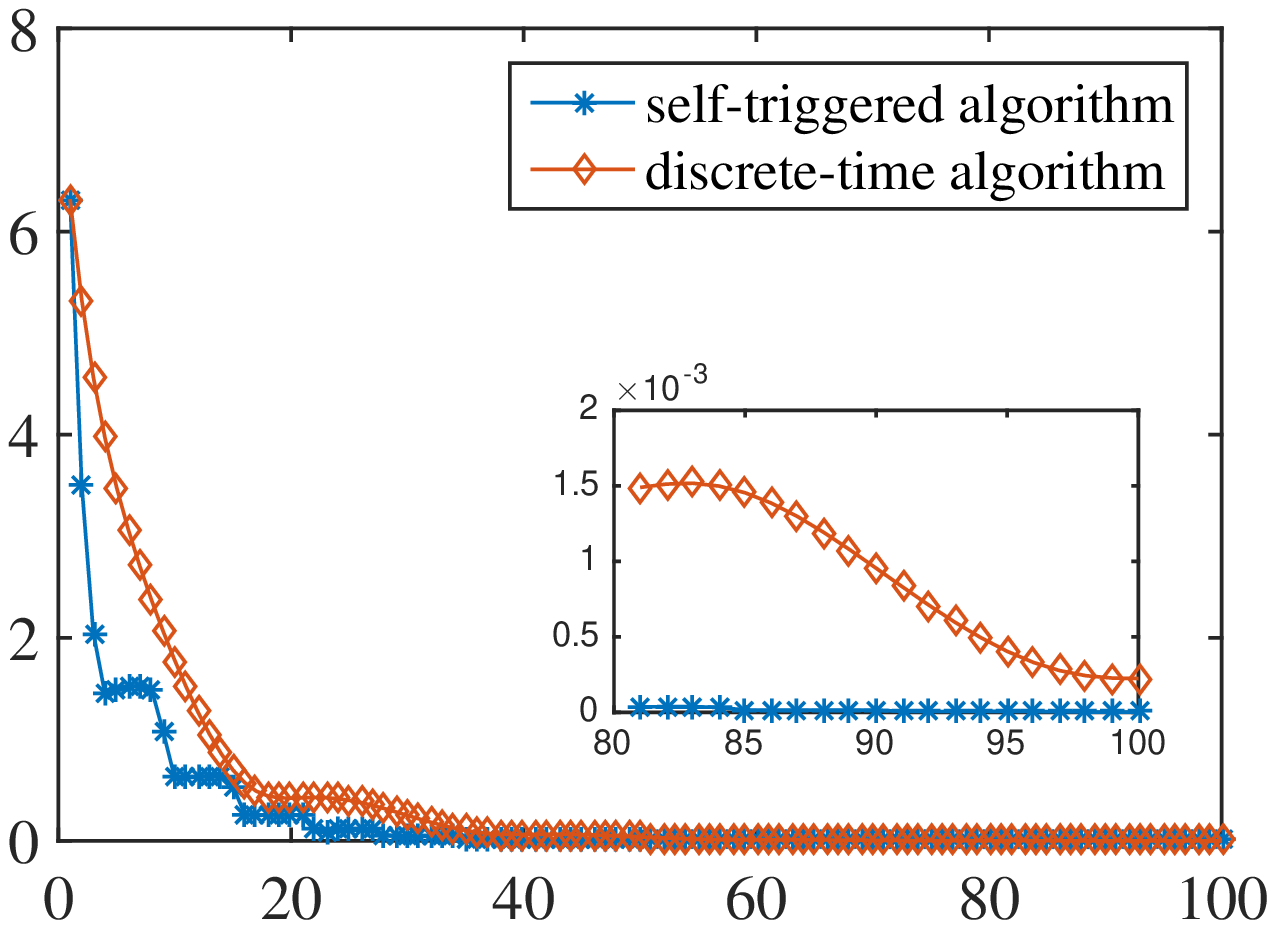}}
    \quad
    \subfloat[$\norm{(x,z)}_{\sdl{F}}$]{\includegraphics[width = 0.475
      \linewidth]{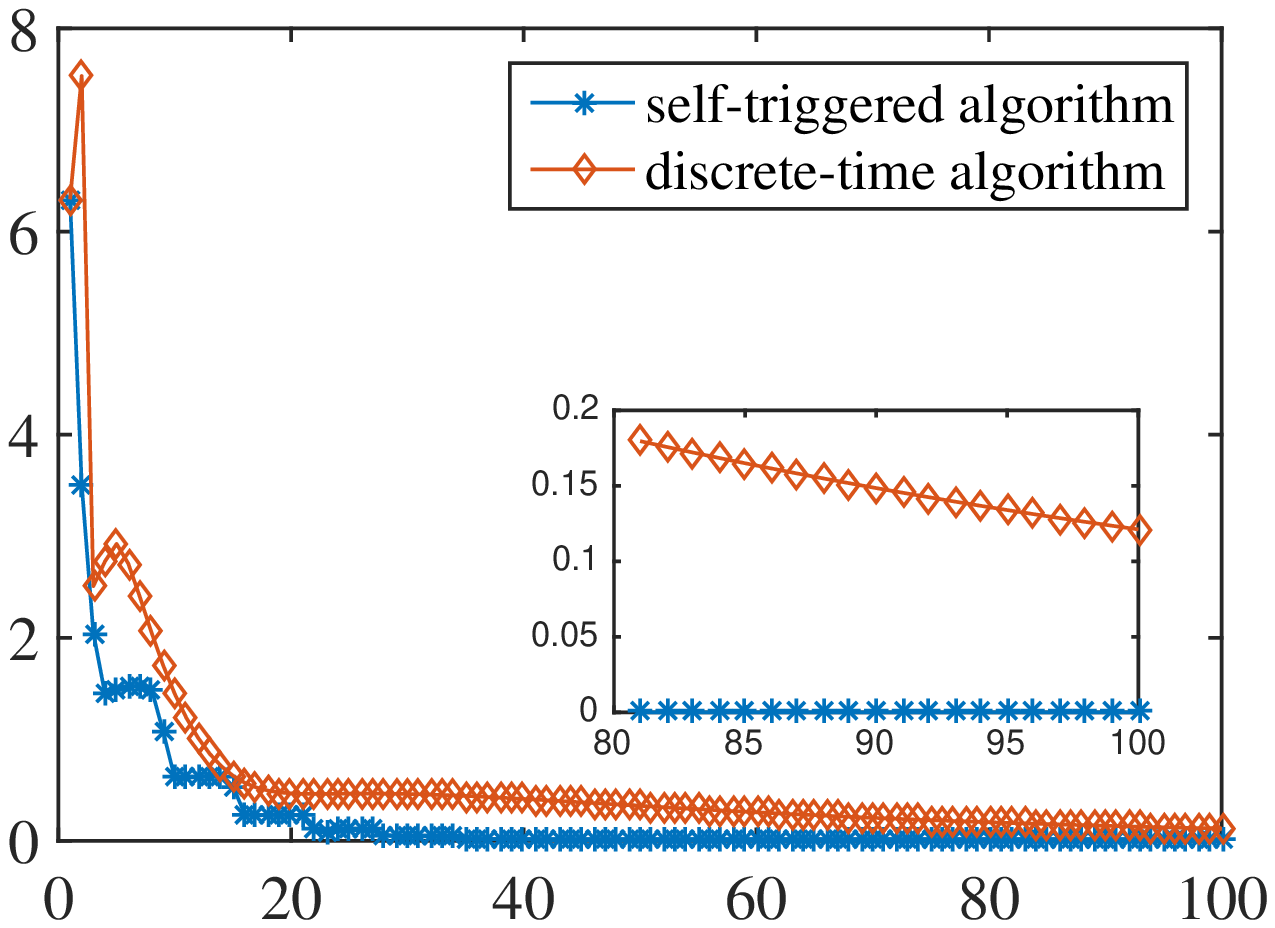}}
    \caption{Comparison between the self-triggered
      saddle-point dynamics and a first-order Euler discretization of
      the saddle-point dynamics with two different stepsize rules. The
      initial condition and implementation details are the same as in
      Figure~\ref{fig:self-triggered}. Both plots show the evolution of
      the distance to the saddle point, compared in (a) against a
      constant-stepsize implementation with value $0.1$ and in (b)
      against a decaying-stepsize implementation with value $1/k$ at the
      $k$-th iteration. The self-triggered dynamics converges faster in
      both cases.}\label{fig:comparison}
  \end{figure}
  \oprocend
}
\end{example}

\section{Conclusions}\label{sec:conclusions}
This paper has studied the global convergence and robustness
properties of the projected saddle-point dynamics.  We have provided a
characterization of the omega-limit set in terms of the Hessian blocks
of the saddle function. Building on this result, we have established
global asymptotic convergence assuming only local strong
convexity-concavity of the saddle function.  For the case when this
strong convexity-concavity property is global, we have identified a
Lyapunov function for the dynamics. In addition, when the saddle
function takes the form of a Lagrangian of an equality constrained
optimization problem, we have established the input-to-state stability
of the saddle-point dynamics by identifying an ISS Lyapunov function,
which we have used to design a self-triggered discrete-time
implementation.  In the future, we aim to generalize the ISS results
to more general classes of saddle functions. In particular, we wish to
define a ``semi-global'' ISS property that we conjecture will hold for
the saddle-point dynamics when we relax the global upper bound on the
Hessian block of the saddle function. Further, to extend the ISS
results to the projected saddle-point dynamics, we plan to develop the
theory of ISS for general projected dynamical systems. Finally, we
intend to apply these theoretical guarantees to determine robustness
margins and design opportunistic state-triggered implementations for
frequency regulation controllers in power networks.

\appendix
\renewcommand{\theequation}{A.\arabic{equation}}
\renewcommand{\thetheorem}{A.\arabic{theorem}}
\section{Appendix}

Here we collect a couple of auxiliary results used in the proof of
Theorem~\ref{th:iss-lyapunov}.

\begin{lemma}\longthmtitle{Auxiliary result for
    Theorem~\ref{th:iss-lyapunov}: I}\label{le:U-bound-new}
  Let $B_1, B_2 \in \real^{n \times n}$ be symmetric matrices satisfying
  $mI \preceq B_1, B_2 \preceq MI$ for some $0 < m \le M < \infty$. Let
  $\beta_1 > 0$, $\beta_2 = \frac{4 \beta_1 M^4}{m^2}$, and $\lm_m =
  \min\set{\frac{1}{2}\beta_1 m, \beta_1 m^3}$.  Then, 
  \begin{align*}
    W:=
    \begin{bmatrix}
      \beta_1 B_1 B_2 B_1 + \beta_2 B_1 & \beta_1 B_1 B_2
      \\
      \beta_1 B_2 B_1 & \beta_1 B_2
      \end{bmatrix}
      \succ \lm_m I.
  \end{align*}
\end{lemma}
\begin{IEEEproof}
  Reasoning with Schur complement~\cite[Section A.5.5]{SB-LV:09}, the
  expression $W - \lm_m I \succ 0$ holds if and only if the following
  hold
  \begin{align}\label{eq:LMI1}
    & \beta_1 B_1 B_2 B_1 + \beta_2 B_1 - \lm_m I\succ 0, \notag
    \\
    & \beta_1 B_2 - \lm_m I -
    \\
    & \, \beta_1 B_2 B_1 (\beta_1 B_1 B_2 B_1
    + \beta_2 B_1 - \lm_m I)^{-1} \beta_1 B_1 B_2 \succ 0.
    \notag
  \end{align}
  The first of the above inequalities is true since $\beta_1 B_1
  B_2 B_1 + \beta_2 B_1 - \lm_m I \succeq \beta_1 m^3 I + \beta_2 m I -
  \lm_m I \succ 0$ as $\lm_m \le \beta_1 m^3$.  For the second
  inequality note that
  \begin{align*}
   & \beta_1 B_2  - \lm_m I \notag
    \\
    & \quad - \beta_1 B_2 B_1 (\beta_1 B_1 B_2 B_1 + \beta_2 B_1
    - \lm_m I)^{-1} \beta_1 B_1 B_2 \notag
    \\ 
    & \succeq (\beta_1 m  - \lm_m) I \notag
    \\
    & \quad - \beta_1^2 M^4 \lambda_{\max}\Bigl((\beta_1
    B_1 B_2 B_1 + \beta_2 B_1 - \lm_m I)^{-1}\Bigr) I \notag
    \\
    & \succeq \Bigl( \frac{1}{2} \beta_1 m  -
    \frac{\beta_1^2  M^4}{\lambda_{\min} (\beta_1 B_1
  B_2 B_1 + \beta_2 B_1 - \lm_m I)} \Bigr) I,
  \end{align*}
  where in the last inequality we have used the fact that $\lm_m \le
  \beta_1 m/2$.  Note that $\lambda_{\min}\Bigl(\beta_1
  B_1 B_2 B_1 + \beta_2 B_1 - \lm_m I\Bigr) \ge \beta_1 m^3 +
  \beta_2 m - \lm_m \ge \beta_2 m$.  Using this lower bound, the
  following holds 
  \begin{align*}
    \frac{1}{2} \beta_1 m - \frac{\beta_1^2 
  M^4}{\lambda_{\min} (\beta_1 B_1 B_2 B_1 + \beta_2 B_1 - \lm_m
  I)}
    \\
    \ge \frac{1}{2} \beta_1 m  - \frac{\beta_1^2
       M^4}{\beta_2 m} = \frac{1}{4} \beta_1 m.
  \end{align*}
  The above set of inequalities show that the second inequality
  in~\eqref{eq:LMI1} holds, which concludes the proof. 
\end{IEEEproof}

\begin{lemma}\longthmtitle{Auxiliary result for
    Theorem~\ref{th:iss-lyapunov}: II}\label{le:A-z-bound}
  Let $F$ be of the form~\eqref{eq:F-wo-y} with $f$ strongly convex. Let
  $(x,z) \in \real^n \times \real^m$ and $(\xo,\zo) =
  \proj_{\sdl{F}}(x,z)$. Then, $z-\zo$ is orthogonal to the kernel of
  $A^\top$, and
  \begin{align*}
    \norm{A^\top (z-\zo)}^2 \ge \lm_{\mathrm{s}}(AA^\top)
    \norm{z-\zo}^2,
  \end{align*}
  where $\lm_{\mathrm{s}}(A A^\top)$ is the smallest nonzero eigenvalue of
  $AA^\top$.
\end{lemma}
\begin{IEEEproof}
  Our first step is to show that there exists $\xo \in \real^n$ such
  that if $(x,z) \in \sdl{F}$, then $x = \xo$. By contradiction, assume
  that $(x_1,z_1), (x_2,z_2) \in \sdl{F}$ and $x_1 \not = x_2$.  The
  saddle point property at $(x_1,z_1)$ and $(x_2,z_2)$ yields
  \begin{align*}
    F(x_1,z_1) \le F(x_2,z_1) \le F(x_2,z_2) 
    \le F(x_1,z_2) \le
    F(x_1,z_1).
  \end{align*}
  This implies that $F(x_1,z_1) = F(x_2,z_1)$, which is a
  contradiction as $x \mapsto F(x,z_1)$ is strongly convex and $x_1$
  is a minimizer of this map. Therefore, $\sdl{F} = \{\xo\} \times
  \ZZ$, $\ZZ \subset \real^m$. Further, recall that the set of saddle
  points of $F$ are the set of equilibrium points of the saddle point
  dynamics~\eqref{eq:saddle-dyn}. Hence, $(\xo,z) \in \sdl{F}$ if and
  only if
  \begin{align*}
    \gradient f(\xo) + A^\top z = 0.
  \end{align*}
  We conclude from this that  
  \begin{align}\label{eq:Z-set}
    \ZZ = -(A^\top)^\dagger \gradient f(\xo) + \ker(A^\top), 
   \end{align}
   where $(A^\top)^\dagger$ and $\ker(A^\top)$ are the Moore-Penrose
   pseudoinverse~\cite[Section A.5.4]{SB-LV:09} and the kernel of
   $A^\top$, respectively. By definition of the projection operator, if
   $(\xo,\zo) = \proj_{\sdl{F}}(x,z)$, then $\zo = \proj_{\ZZ}(z)$ and
   so, from~\eqref{eq:Z-set}, we deduce that $(z-\zo)^\top v = 0$ for
   all $v \in \ker(A^\top)$. Using this fact, we conclude the proof by
   writing
   \begin{align*}
     \norm{A^\top (z-\zo)}^2 = (z-\zo)^\top A A^\top (z - \zo) 
     \\
     \ge
     \lm_{\mathrm{s}}(A A^\top) \norm{z - \zo}^2,
   \end{align*}
   where the inequality follows by writing the eigenvalue
   decomposition of $AA^\top$, expanding the quadratic expression in
   $(z-\zo)$, and lower-bounding the terms. 
\end{IEEEproof}

\begin{proposition}\longthmtitle{Gradient of $V_3$ is locally
    Lipschitz}\label{prop:V-bound}
  Let the saddle function $F$ be of the form~\eqref{eq:F-wo-y}, with $f$
  twice differentiable, map $x \mapsto \gradient^2 f(x)$ Lipschitz
  with some constant $L>0$, and $m I \preceq \gradient^2 f(x) \preceq M
  I$ for all $x \in \real^n$ and some constants $0< m \le M < \infty$.
  Then, for $V_3$ given in~\eqref{eq:V2}, the following holds
  \begin{align*}
    \norm{\gradient V_3(x_2,z_2) - \gradient V_3(x_1,z_1)} \le
    \xi(x_1,z_1) \norm{x_2-x_1; z_2-z_1},
  \end{align*}
  for all $(x_1,z_1), (x_2,z_2) \in \real^n \times \real^m$, where 
  \begin{align}
    \xi (x_1,z_1) & = \sqrt{3} \Bigl( \beta_1^2 (\xi_1(x_1,z_1)^2 + \norm{A}^4
    + \norm{A}^2 \xi_2^2) + \beta_2^2 \Bigr)^{\frac{1}{2}}, \notag
    \\
    \xi_1(x_1,z_1) & = M \xi_2 + L \norm{\gradient_x F(x_1,z_1)},
    \notag
    \\
    \xi_2 & = \max\set{M,\norm{A}}. \label{eq:xi-def}
  \end{align}
\end{proposition}
\begin{IEEEproof}
  For the map $(x,z) \mapsto \gradient_x F(x,z)$, note that
  \begin{align}
    & \norm{\gradient_x F(x_2,z_2) - \gradient_x F(x_1,z_1)} \notag
    \\ 
    & \quad = \Bnorm{ \int_0^1 \gradient_{xx} F(x(s),z(s)) (x_2 - x_1)
    ds \notag
    \\
    & \qquad \qquad \qquad \qquad \qquad  + \int_0^1 \gradient_{zx}
  F(x(s),z(s)) (z_2 - z_1) } \notag
    \\
    & \quad \le M \norm{x_2 - x_1} + \norm{A} \norm{z_2 - z_1} \notag
    \\
    & \quad \le \xi_2 \norm{x_2-x_1;z_2-z_1}, \label{eq:lip-c1}
  \end{align}
  where $x(s) = x_1 + s(x_2 - x_1)$, $z(s) = z_1 + s(z_2 - z_1)$ and
  $\xi_2 = \max\set{M,\norm{A}}$.
  In the above inequalities we have
  used the fact that $\norm{\gradient_{xx} F(x,z)} = \norm{\gradient^2
  f(x)} \le M$ for any $(x,z)$. Further, the following Lipschitz
  condition holds by assumption
  \begin{align}
    \norm{\gradient_{xx} F(x_2,z_2) - \gradient_{xx} F(x_1,z_1)} \le L
    \norm{x_2 - x_1} \label{eq:lip-c2}
  \end{align}
  Using~\eqref{eq:lip-c1} and~\eqref{eq:lip-c2}, we get
  \begin{align}
    & \norm{\gradient_{xx} F(x_2,z_2)\gradient_xF(x_2,z_2) -
    \gradient_{xx} F(x_1,z_1) \gradient_x F(x_1,z_1)}  \notag
    \\
    & \quad \le \norm{\gradient_{xx}F(x_2,z_2) (\gradient_x F(x_2,z_2)
    - \gradient_x F(x_1,z_1))} \notag
    \\ 
    & \quad \quad + \norm{(\gradient_{xx} F(x_2,z_2) - \gradient_{xx}
    F(x_1,z_1)) \gradient_x F(x_1,z_1)} \notag
    \\
    & \quad \le \xi_1(x_1,z_1) \norm{x_2 - x_1; z_2 - z_1},
    \label{eq:lip-c3}
  \end{align}
  where $\xi_1(x_1,z_1) = M \xi_2 + L \norm{\gradient_x F(x_1,z_1)}$.
  Also, 
  \begin{align}
     \norm{\gradient_z F(x_2,z_2)  - & \gradient_z F(x_1,z_1)} 
      = \norm{A(x_2 - x_1)} \notag
     \\
     & \le \norm{A} \norm{x_2 - x_1;z_2 - z_1} \label{eq:lip-c4}
  \end{align}
  Now note that 
  \begin{align*}
    \gradient_x V_3(x,z) & = \beta_1 \Bigl(\gradient_{xx} F(x,z)
    \gradient_x F(x,z) + A^\top \gradient_z F(x,z) \Bigr) 
    \\
    &  \qquad + \beta_2 (x -\xo),
    \\
    \gradient_z V_3(x,z) & = \beta_1 A \gradient_x F(x,z) + \beta_2 (z
    - \zo).&
  \end{align*}
  Finally, using~\eqref{eq:lip-c1},~\eqref{eq:lip-c3},
  and~\eqref{eq:lip-c4}, we get
  \begin{align*}
    & \norm{\gradient V_3(x_2,z_2) - \gradient V_3(x_1,z_1)}^2 =
    \norm{\gradient_x V_3(x_2,z_2) 
    \\ 
    & \qquad - \gradient_x V_3(x_1,z_1)}^2 + \norm{\gradient_z
    V_3(x_2,z_2) - \gradient_z V_3(x_1,z_1)}^2
    \\
    & \overset{(a)}{\le} 3 \beta_1^2 \norm{\gradient_{xx} F(x_2,z_2) \gradient_x
    F(x_2,z_2) 
    \\
    & \qquad \qquad \qquad \qquad \qquad \quad - \gradient_{xx}
    F(x_1,z_1) \gradient_x F(x_1,z_1)}^2 
    \\
    & \,  + 3 \beta_1^2 \norm{A^\top(\gradient_z F(x_2,z_2) -
    \gradient_z F(x_1,z_1))}^2 \! + \! 3 \beta_2^2 \norm{x_2 - x_1}^2
    \\
    & \,  + 3 \beta_1^2 \norm{A (\gradient_x F(x_2,z_2) - \gradient_x
    F(x_1,z_1))}^2 \! + \!3 \beta_2^2 \norm{z_2 - z_1}^2
    \\
    & \le \xi(x_1,z_1)^2 \norm{x_2 - x_1; z_2-z_1}^2,
  \end{align*}
  where in (a), we have used the inequality $(a+b)^2 \le 3(a^2+b^2)$
  for any $a,b \in \real$. This concludes the proof.  
\end{IEEEproof}

\section*{Acknowledgments}
We would like to thank Simon K. Niederl\"ander for discussions on
Lyapunov functions for the saddle-point dynamics and \changes{an
  anonymous reviewer for suggesting the ISS Lyapunov function of
  Remark~\ref{re:alternative-ISS}}.

\bibliographystyle{IEEEtran}
\bibliography{alias,FB,JC,Main,Main-add}

\end{document}